\documentclass[12pt,a4paper,reqno]{amsart}
\usepackage{amssymb}
\usepackage{latexsym}
\usepackage{exscale}
\usepackage{mathrsfs}

\newcommand{\cB}{\mathcal B}

\newcommand{\be}{{\mathbf e}}
\newcommand{\bs}{{\mathbf s}}
\newcommand{\cD}{{\mathcal{D}}}

\newcommand {\SB} {{\mathbb B}}

\newcommand {\SE} {{\mathbb E}}

\newcommand {\SN} {{\mathbb N}}
\newcommand {\SR} {{\mathbb R}}
\newcommand {\SW} {{\mathbb W}}
\newcommand {\SX} {{\mathbb X}}

\newcommand {\al} {{\alpha}}
\newcommand {\dt} {{\delta}}
\newcommand {\Dt} {{\Delta}}

\newcommand {\ga} {{\gamma}}
\newcommand {\Ga} {{\Gamma}}
\newcommand {\la} {{\lambda}}

\newcommand {\om} {{\omega}}
\newcommand{\cA}{\mathcal A}

\numberwithin{equation}{section}
\newtheorem{theorem}{Theorem}[section]
\newtheorem{lemma}[theorem]{Lemma}

\newtheorem{corollary}[theorem]{Corollary}
\newtheorem{Remark}[theorem]{Remark}
\newtheorem{proposition}[theorem]{Proposition}

\newtheorem{example}[theorem]{Example}

\newcommand{\Ba}[1]{\begin{array}{#1}}
\newcommand{\Ea}{\end{array}}
\newcommand{\Be}{\begin{equation}}
\newcommand{\Ee}{\end{equation}}
\newcommand{\Bea}{\begin{eqnarray}}
\newcommand{\Eea}{\end{eqnarray}}
\newcommand{\Beas}{\begin{eqnarray*}}
\newcommand{\Eeas}{\end{eqnarray*}}
\newcommand{\Benu}{\begin{enumerate}}
\newcommand{\Eenu}{\end{enumerate}}

\newcommand{\BR}{\begin{Remark} \em}
\newcommand{\ER}{\end{Remark}}
\newcommand{\BE}{\begin{example} \em}
\newcommand{\EE}{\end{example}}

\newcommand {\Ds} {\displaystyle}
\newcommand {\Ts} {\textstyle}

\newcommand {\mand} {{\quad\mbox{and}\quad}}
\renewcommand {\mid} {{\,\,\,\colon\,\,\,}}

\newcommand{\bline}{{\bigskip

\noindent}}

\newcommand{\sline}{{\smallskip

\noindent}}

\newcommand{\fpr}{{{\mathfrak f}^s_{p,r}}}
\newcommand{\hpr}{{{\mathfrak f}^p_{\rm{hyp}} }}
\newcommand{\lpr}{{{\mathfrak l}^{p,q}}}

\newcommand{\fphi}{{{\mathfrak f}^{\Phi}}}

\newcommand{\cW}{{\mathcal{W}}}

\newcommand {\bone} {{\bf 1}}
\newcommand{\Gaq}{{\mathscr{G}_{q}^{\alpha}}}
\newcommand{\GaqB}{{\mathscr{G}_{q}^{\alpha}(\mathcal{B},\mathbb{B})}}
\newcommand{\Aaq}{{\mathcal{A}_{q}^{\alpha}}}
\newcommand{\AaqB}{{\mathcal{A}_{q}^{\alpha}(\mathcal{B},\mathbb{B})}}

\begin{document}
\allowdisplaybreaks

\title[Democracy and embeddings]{Democracy Functions and Optimal Embeddings  for
Approximation spaces}

\author{Gustavo Garrig\'os}

\address{Gustavo Garrig\'os
\\
Departamento de Matem\'aticas
\\
Universidad Aut\'onoma de Ma\-drid
\\
28049, Madrid, Spain}

\email{gustavo.garrigos@uam.es}

\author{Eugenio Hern\'andez}

\address{Eugenio, Hern\'andez
\\
Departamento de Matem\'aticas
\\
Universidad Aut\'onoma de Madrid
\\
28049, Madrid, Spain}

\email{eugenio.hernandez@uam.es}

\author{Maria de Natividade}

\address{Maria de Natividade
\\
Departamento de Matem\'aticas
\\
Universidad Aut\'onoma de Madrid
\\
28049, Madrid, Spain}

\email{maria.denatividade@uam.es}

\begin{abstract}
We prove optimal embeddings for nonlinear approximation spaces
$\mathcal{A}^\al_q$, in terms of weighted Lorentz sequence spaces,
with the weights depending on the democracy functions of the
basis. As applications we recover known embeddings for $N$-term
wavelet approximation in $L^p$, Orlicz, and Lorentz norms. We also
study the ``greedy classes'' $\Gaq$ introduced by Gribonval and
Nielsen,   obtaining new counterexamples which show that
$\Gaq\not=\mathcal{A}^\al_q$ for most non democratic unconditional
bases.
\end{abstract}

\thanks{Research supported by Grant MTM2007-60952 of Spain. The research
of M. de Natividade supported by Instituto Nacional de Bolsas de
Estudos de Angola, INABE}

\date{\today}
\subjclass[2000]{41A17, 42C40.}

\keywords{Non-linear approximation, greedy algorithm, democratic
bases, Jackson and Bernstein inequalities, discrete Lorentz
spaces, wavelets. }

\maketitle

\section{Introduction}\label{secIntroduc}

Let $(\mathbb{B},\|.\|_{\mathbb{B}})$ be a quasi-Banach space with
a countable {\bf unconditional} basis  $\mathcal{B} = \{e_{j}: j
\in \mathbb{N}\}.$ A main question in {\bf Approximation Theory}
consists in finding a characterization (if possible) or at least
suitable embeddings for the non-linear approximation spaces
$\mathcal{A}_{q}^{\alpha}(\mathcal{B},\mathbb{B})$, $\al>0$,
$0<q\leq\infty$, defined using the {\bf N-term error of
approximation} $\sigma_N(x, \mathbb{B})$ (see sections
\ref{subNonLinAppGreedyAlgor} and \ref{sbsApSGredS} for
definitions). Such characterizations or inclusions are often given
in terms of ``smoothness classes'' of the sort
\[
\mathfrak{b}(\cB;\SB)\,:=\,
\left\{x=\sum_{j=1}^{\infty} c_j e_j \in\SB\;\colon\; \{\|c_je_j\|_\SB\}_{j=1}^\infty\in\mathfrak{b}\right\},
\]
where $\mathfrak{b}$ is a suitable sequence space whose elements
decay at infinity, such as $\ell^\tau$ or more generally the
discrete Lorentz classes $\ell^{\tau,q}$.

The simplest result in this direction appears  when $\mathcal{B}$
is an orthonormal basis in a Hilbert space $\mathbb{H}\,$, and was
first proved by Stechkin when $\al=1/2$ and $q=1$ (see \cite{ST}
or \cite{DT} for general $\al,q$).

\begin{theorem}\label{StDT} {\rm (\cite{ST,DT}).}
  Let $\mathcal{B} = \{e_j\}_{j=1}^\infty$ be an orthonormal basis
in a Hilbert space $\mathbb{H}$, and $\alpha > 0$, $ 0 < q
\leq \infty$. Then
 $$ \mathcal{A}_{q}^{\alpha}(\mathcal{B},\mathbb{H}) =
 \ell^{\tau, q}(\mathcal{B};\mathbb{H})$$
where $\tau$ is defined by $\frac{1}{\tau} = \alpha + \frac12$.
\end{theorem}


Many results have been published in the literature similar to
Theorem \ref{StDT} when  $\mathbb H$ is replaced by
a particular space (say, $L^p$) and the basis $\mathcal B$
is a particular one (for example, a wavelet basis). We refer to the
survey articles \cite{DeV} and \cite{Tem} for detailed statements and
references.

There are also a number of results for general pairs $(\SB,\cB)$
(even with the weaker notion of quasi-greedy basis
\cite{GN,DKKT,KP}). We recall two of them in the setting of
unconditional bases which we consider here. For simplicity, in all
the statements we assume that the basis is \emph{normalized},
meaning $\|e_j\|_\SB=1$, $\forall\,j\in\SN$. The first result can
be found in \cite{KP2} (see also \cite{GH}).

\begin{theorem}\label{KPGH} {\rm (\cite[Th 1]{KP2}, \cite[Th 6.1]{GH}).}
 Let $\mathbb{B}$ be a quasi-Banach space
and $\mathcal B = \{e_{j}\}_{j=1}^\infty$ a (normalized)
unconditional basis satisfying the following property: there exists
$p \in (0, \infty)$ and a constant $C > 0$ such that
  \begin{eqnarray} \label{p-space}
   \frac{1}{C}|\Gamma|^{1/p}\leq \Big\| \sum_{k\in \Gamma}
   e_k \Big\|_{\mathbb B} \leq C |\Gamma|^{1/p}
  \end{eqnarray}
for all finite $\Gamma \subset \mathbb N\,.$ Then, for $\al>0$ and $ 0 < q \leq \infty$ we have
  \[\mathcal{A}_{q}^{\alpha}(\mathcal{B},\mathbb{B}) =
 \ell^{\tau, q}(\mathcal{B};\mathbb{B})
\]
when $\tau$ is defined by $\frac{1}{\tau} = \alpha + \frac1p$.
\end{theorem}

Condition (\ref{p-space}) is sometimes referred as $\cB$ having
the $p$-Temlyakov property \cite{KP}, or as $\mathbb B$ being a
$p$-space \cite{HJLY,GH}. For instance, wavelet bases in $L^p$
satisfy this property \cite{Tem1}. The second result we quote is proved in
\cite{GN} (see also \cite{KP2}).

\begin{theorem}\label{GN} {\rm (\cite[Th 3.1]{GN}).}
 Let $\mathbb{B}$ be a Banach space
and $\mathcal B = \{e_{j}\}_{j=1}^\infty$ a (normalized)
unconditional basis with the following property: there exist $1
\leq p \leq q \leq \infty$ and constants $A,B>0$ such that when $x
= \sum_{j\in \mathbb N} c_j e_j \in \mathbb B$ we have
  \begin{eqnarray}\label{sandwich}
  A \,\| \{c_j\}\|_{\ell^{q,\infty}}\, \leq \,\| x \|_{\mathbb B} \,\leq B \,\|
  \{c_j\}\|_{\ell^{p,1}}\,.
  \end{eqnarray}
Then, for $\al>0$ and $ 0 < s \leq \infty$ we have
  \begin{eqnarray}\label{inclusions}
  \ell^{\tau_p, s}(\mathcal{B};\mathbb{B}) \hookrightarrow
  \mathcal{A}_{s}^{\alpha}(\mathcal{B},\mathbb{B})
 \hookrightarrow \ell^{\tau_q, s}(\mathcal{B};\mathbb{B})
 \end{eqnarray}
where $\frac{1}{\tau_p} = \alpha + \frac1p $ and
$\frac{1}{\tau_q} = \alpha + \frac1q $. Moreover, the inclusions
given in (\ref{inclusions}) are best possible in the sense described
in section 4 of \cite{GN}.
\end{theorem}

Condition (\ref{sandwich}) is referred in \cite{GN} as  $(\mathbb
B , \mathcal B)$ having the $(p,q)$ {\bf sandwich property}, and
it is shown to be equivalent to
  \begin{eqnarray}\label{pq-space}
  A' |\Gamma|^{1/q} \leq \Big\| \sum_{k\in \Gamma}
 e_k\Big\|_{\mathbb B} \leq B' |\Gamma|^{1/p}
  \end{eqnarray}
for all $\Ga\subset\SN$ finite. Observe that \eqref{pq-space}
coincides with (\ref{p-space}) when $p = q\,.$

The purpose of this article is to obtain optimal embeddings for $\cA^\al_q(\cB,\SB)$
as in \eqref{inclusions} when no condition such as \eqref{pq-space}
is imposed. More precisely, we define
the {\bf right} and {\bf left democracy functions} associated with a 
basis $\mathcal B$ in $\SB$ by
 $$h_{r}(N;\mathcal{B},\mathbb{B})\equiv
\sup_{|\Gamma|=N}\Big\|\sum_{k\in\Gamma}\frac{e_{k}}{\|e_k\|_\SB}\Big\|_{\mathbb{B}}
\quad \mbox{and} \quad
 h_\ell(N;\mathcal{B},\mathbb{B})\equiv \inf_{|\Gamma| =
N}\Big\|\sum_{k\in\Gamma}\frac{e_{k}}{\|e_k\|_\SB}\Big\|_{\mathbb{B}}
  $$
for $N=1,2,3, \dots$. We refer to section \ref{Examples} for
various examples where $h_\ell(N)$ and $h_r(N)$ are computed
explicitly (modulo multiplicative constants). As usual, when
$h_\ell(N)\approx h_r(N)$ for all $N\in\SN$ we say that $\cB$ is a
democratic basis in $\SB$ \cite{KT}.

The embeddings will be given in terms of  {\bf weighted discrete
Lorentz spaces} {\boldmath$\ell_{\eta}^{\,q}$}, with quasi-norms
defined by
\[
\bigl\|\{c_k\}\bigr\|_{\ell^{\,q}_\eta}\,\equiv\,
\Big(\,\sum_{k=1}^{\infty}
\,\bigl|\eta(k)\,c^*_k\bigr|^{q}\,\tfrac{1}{k}\,\Big)^{\frac{1}{q}},
\]
where $\{c^*_k\}$ denotes the decreasing rearrangement of
$\{|c_k|\}$ and the \emph{weight} $\eta =
\{\eta(k)\}_{k=1}^\infty$ is a suitable sequence
increasing to infinity and satisfying the doubling property (see
section \ref{subDiscLorSpac} for precise definitions and
references). In the special case $\eta(k) = k^{1/\tau}$ we
 recover the classical definition $\ell_{\eta}^q  = \ell^{\tau,q}$.

\begin{theorem}\label{InclusionsTheorem}
 Let $\mathbb{B}$ be a quasi-Banach space
and $\mathcal B$ an unconditional basis. Assume that $h_\ell(N)$
is doubling. Then if $\al>0$ and $ 0 < q \leq \infty$ we have the
continuous embeddings
  \begin{eqnarray}\label{inclusions-general}
  \ell_{k^\alpha h_r(k)}^{q}(\mathcal{B};\mathbb{B}) \hookrightarrow \mathcal{A}_{q}^{\alpha}(\mathcal{B},\mathbb{B})
 \hookrightarrow \ell_{k^\alpha
 h_\ell(k)}^{q}(\mathcal{B};\mathbb{B})\,.
 \end{eqnarray}
Moreover, for fixed $\al$ and $q$ these inclusions are best
possible in the scale of weighted discrete Lorentz spaces
$\ell_{\eta}^q $, in the sense explained in sections 3, 4 and 6.
\end{theorem}


Observe that this theorem generalizes Theorems \ref{KPGH} and
\ref{GN}. In Theorem \ref{KPGH} we have $h_r(N) \approx h_\ell(N)
\approx N^{1/p}$ and in Theorem \ref{GN}, $h_r(N) \lesssim N^{1/p}$
and $h_\ell(N) \gtrsim N^{1/q}$. When $\mathcal B$ is democratic in
$\mathbb B$, Theorem \ref{InclusionsTheorem} shows that
$\mathcal{A}_{q}^{\alpha}(\mathcal{B},\mathbb{B}) \approx
\ell_{k^\alpha h(k)}^{q}(\mathcal{B};\mathbb{B})$ with $h(k) =
h_r(k) \approx h_\ell(k)\,.$ Compare this result with Corollary 1 in
\cite[\S6]{GN}.

Theorem \ref{InclusionsTheorem} is a consequence of the results
proved in sections \ref{secRigtDemJackIneq} and
\ref{secLeftDemBernsIneq}. Section \ref{secRigtDemJackIneq} deals
with the lower embedding in (\ref{inclusions-general}) and shows
the relation to Jackson type inequalities. Section
\ref{secLeftDemBernsIneq} deals with the upper embedding of
(\ref{inclusions-general}) and its relation to Bernstein type
inequalities. Section \ref{Examples} contains various examples of
democracy functions and embeddings with precise references; these
are all special cases of Theorem \ref{InclusionsTheorem}. In
section 6 we apply Theorem \ref{InclusionsTheorem} to estimate the
democracy functions $h_\ell$ and $h_r$ of the approximation space
$\cA^\al_q$.

Finally, the last section of the paper is dedicated to study the
\lq\lq greedy classes\rq\rq\ $\Gaq(\mathcal B, \mathbb B)$
introduced by Gribonval and Nielsen in \cite{GN}, and their
relations with the approximation spaces $\cA^\al_q(\cB,\SB)$. The
classes $\Gaq$ are defined similarly to the approximation spaces,
but with the error of approximation $\sigma_N(x)$ replaced by the
quantity $\|x-G_N(x)\|_\SB$ (see section \ref{sbsApSGredS} for
details).
 It is easy to see that $ \Gaq(\mathcal B,
\mathbb B) \subset \mathcal A _q^\alpha(\mathcal B, \mathbb B)$,
and when $\mathcal B$ is democratic, $ \Gaq(\mathcal B, \mathbb B)
= \mathcal A _q^\alpha(\mathcal B, \mathbb B)\,.$ One may
conjecture that for unconditional bases $\cB$ the converse is
true, that is $ \Gaq(\mathcal B, \mathbb B) = \mathcal A
_q^\alpha(\mathcal B, \mathbb B)$ implies $\cB$ democratic. We do
not know how to show this, but we can exhibit a fairly general
class of non democratic pairs $(\mathcal B,\SB)$ for which
$\Gaq(\mathcal B, \mathbb B) \neq \mathcal A _q^\alpha(\mathcal B,
\mathbb B)$  for all $\alpha > 0$ and  $q \in (0, \infty]\,.$
These include wavelet bases in the non democratic settings of
 $L^{p,q}$ and  $L^p (\log
L)^\alpha$. We also illustrate how irregular the classes
$\Gaq(\mathcal B, \mathbb B)$ can be when $\cB$ is not democratic,
showing in simple situations that they are not even linear spaces.

\section{General Setting}\label{SecgenSettings}

\subsection{Bases}\label{subBases}

Since we work in the setting of quasi-Banach spaces
$(\mathbb{B},\|\cdot\|_{\mathbb{B}})$, we shall often use the
$\rho$-power triangle
inequality\begin{eqnarray}\label{rhotrianIneq}\|x+y\|_{\mathbb{B}}^{\rho}\leq
\|x\|_{\mathbb{B}}^{\rho} +
\|y\|_{\mathbb{B}}^{\rho}\,,\end{eqnarray} which holds for a
sufficiently small $\rho=\rho_\SB\in(0,1]$ (and hence for all
$\mu\leq\rho_\SB$); see \cite[Lemma 3.10.1]{BL}. The case
$\rho_\SB= 1$ gives a Banach space.

 A sequence of vectors $\mathcal{B} =
\{e_{j}\}_{j=1}^\infty$ is a basis of $\mathbb{B}$ if every
$x\in\mathbb{B}$ can be uniquely represented as $x =
\sum_{j=1}^\infty c_{j}e_{j}$ for some scalars $c_{j}$, with
convergence in $\|\cdot\|_\SB$. The basis $\cB$ is {\bf
unconditional} if the series converges unconditionally, or
equivalently if there is some $K > 0$ such that
\begin{eqnarray}\label{uncondProper}
\Big\|\sum_{j=1}^\infty\lambda_{j}c_{j}e_{j}\Big\|_{\mathbb{B}}\leq
K \Big\|\sum_{j=1}^\infty c_{j}e_{j}\Big\|_{\mathbb{B}}
\end{eqnarray}
for every sequence of scalars $\{\lambda_j\}_{j=1}^\infty$ with
$|\lambda_j| \leq 1$ (see eg \cite[Chapter 5]{HW}).

For simplicity in the statements, throughout the paper we shall
assume that $\cB$ is a {\bf normalized} basis, meaning
$\|e_j\|_{\mathbb B}=1$ for all $j\in \mathbb N\,.$ We can also
assume that the unconditionality constant in \eqref{uncondProper}
is $K=1$. To see so, one can introduce an equivalent quasi-norm in
$\SB$
 $$
 |\!|\!|x|\!|\!|_{\mathbb B} = \sup_{\Gamma \text{finite}, |\lambda_j|\leq
 1} \|\sum_{j\in \Gamma} \lambda_j x_j e_j\|_{\mathbb B}\,, \qquad \mbox{ if } x
 = \sum_{j=1}^\infty x_j e_j.
 $$
 Observe that with this renorming we still have
$|\!|\!|e_j|\!|\!|_{\mathbb B}=1$.

With the above assumptions, the following \textbf{lattice
property} holds: if $|y_{k}|\leq |x_{k}|$ for all $k\in
\mathbb{N}$ and $x = \sum_{k=1}^{\infty}x_{k}e_{k} \in
\mathbb{B},$ then the series $y = \sum_{k=1}^{\infty}y_{k}e_{k}$ converges in $\mathbb{B}$ and $\|y\|_{\mathbb{B}}\leq \|x\|_{\mathbb{B}}.$ Also, using \eqref{uncondProper} with $K=1$ we see that,
 for every $\Gamma \subset \mathbb{N}$ finite
\begin{eqnarray}\label{LaticeIneq}(\inf_{j\in\Gamma}|c_{j}|)\Big\|\sum_{j\in\Gamma}e_{j}\Big\|_{\mathbb{B}}\leq
\Big\|\sum_{j\in\Gamma}c_{j}e_{j}\Big\|_{\mathbb{B}}\leq
(\sup_{j\in\Gamma}|c_{j}|)\Big\|\sum_{j\in\Gamma}e_{j}\Big\|_{\mathbb{B}}.
\end{eqnarray}

\subsection{Non-Linear Approximation and Greedy
Algorithm}\label{subNonLinAppGreedyAlgor}

Let $\mathcal{B} = \{e_{j}\}_{j=1}^\infty$ be a basis in
$\mathbb{B}$. Let $\Sigma_{N}$, $N=1,2,3,\ldots$, be the set of all
$y\in\mathbb{B}$ with at most $N$ non-null coefficients in the
unique basis representation.
For $x \in\mathbb{B},$ the \textbf{$N$-term error of
approximation} with respect to $\mathcal{B}$ is defined as
\[\sigma_N(x)=\sigma_{N}(x;\mathcal{B},\mathbb{B})\equiv
\inf_{y\in\Sigma_{N}}\|x-y\|_{\mathbb{B}},\quad N=1,2,3\ldots
\]
We also set  $\Sigma_0=\{0\}$ so that $\sigma_0 (x) =
\|x\|_{\mathbb B} \,.$ Using the lattice property mentioned in
$\S$2.1 it is easy to see that for $x =
\sum_{j=1}^{\infty}c_{j}e_{j}$ we actually have
\begin{eqnarray}\label{ErrorNterApIneq}\sigma_{N}(x)\,=\,
\inf_{|\Gamma| =
N}\Big\{\big\|x-\sum_{\ga\in\Gamma}c_{\ga}e_{\ga}\big\|_{\mathbb{B}}\Big\},
\end{eqnarray}
that is, only coefficients from $x$ are relevant when computing
$\sigma_N(x)$; see eg \cite[(2.6)]{GH}.

 Given
$x=\sum_{j=1}^{\infty}c_{j}e_{j} \in \mathbb{B},$ let $\pi$ denote
any bijection of $\mathbb{N}$ such that
\begin{eqnarray}\label{NonincreRearrang}\|c_{\pi(j)}e_{\pi(j)}\|\geq
\|c_{\pi(j+1)}e_{\pi(j+1)}\|, \quad \mbox{for all} \quad
j\in\mathbb{N}.
\end{eqnarray}
Without loss of generality we may assume that the basis is
normalized and then (\ref{NonincreRearrang}) becames $
|c_{\pi(j)}|\geq |c_{\pi(j+1)}|, \quad \mbox{for all} \quad
j\in\mathbb{N}.$
 A \textbf{greedy algorithm of step $N$} is
a correspondence assigning
\[x =
\sum_{j=1}^{\infty}c_{j}e_{j}\in \mathbb{B}\longmapsto
G_{N}^{\pi}(x) \equiv
\sum_{j=1}^{N}c_{\pi(j)}e_{\pi(j)}
\]
for any $\pi$ as in \eqref{NonincreRearrang}. The
\textbf{error of greedy approximation} at step $N$ is defined by
\begin{eqnarray}\label{defnerrorAppGreedy}\ga_N(x)=\gamma_{N}(x;\mathcal{B},\mathbb{B})
\equiv\sup_{\pi}\|x-G_{N}^{\pi}
(x)\|_{\mathbb{B}}.
\end{eqnarray}
Notice that $\sigma_{N}(x)\leq \ga_N(x)$, but the reverse
inequality may not be true in general. It is said that
$\mathcal{B}$ is a \textbf{greedy basis} in $\mathbb{B}$ when
there is a constant $c \geq 1$ such that
\[\ga_N(x;\mathcal{B},\mathbb{B})\,\leq\,c\,
\sigma_{N}(x;\mathcal{B},\mathbb{B}),\quad\forall\;x\in\SB,\;N=1,2,3,\ldots
\]
A celebrated theorem of Konyagin and Temlyakov characterizes
greedy bases as those which are unconditional and democratic \cite{KT}.

\subsection{Approximation Spaces and Greedy
Classes}\label{sbsApSGredS}

The classical non-linear approximation spaces
$\mathcal{A}_{q}^{\alpha}(\mathcal{B},\mathbb{B})$ are defined as
follows: for $\alpha
> 0$ and $0<q<\infty$
\[\mathcal{A}_{q}^{\alpha}(\mathcal{B},\mathbb{B}) =
\Big\{x\in\mathbb{B}:\|x\|_{\mathcal{A}^{\alpha}_{q}}\equiv
\|x\|_{\mathbb{B}} +
\Big[\sum_{n=1}^{\infty}\big(N^{\alpha}\sigma_{N}(x;\mathcal{B},\mathbb{B})\big)^{q}\frac{1}{N}\Big]^{\frac{1}{q}}<\infty\Big\}.
\]
When $q=\infty$ the  definition takes the form:
\[\mathcal{A}_{\infty}^{\alpha}(\mathcal{B},\mathbb{B}) =
\big\{x\in\mathbb{B}:
\|x\|_{\mathcal{A}^{\alpha}_{\infty}}\equiv\|x\|_{\mathbb{B}} +
\sup_{N\geq
1}N^{\alpha}\sigma_{N}(x)<\infty\big\}.
\]

It is well known that
$\mathcal{A}_{q}^{\alpha}(\mathcal{B},\mathbb{B})$  are
quasi-Banach spaces (see eg \cite{Pie}). Also, equivalent
quasi-norms can be obtained restricting to dyadic $N$'s:
\[\|x\|_{\mathcal{A}_{q}^{\alpha}}\approx \|x\|_{\mathbb{B}} +
\Big[\sum_{k=0}^{\infty}\bigl(2^{k\alpha}\sigma_{2^{k}}(x)\bigr)^{q}\Big]^{\frac{1}{q}}
\]
and likewise for $q=\infty$.
This is a simple consequence of the monotonicity of $\sigma_N(x)$
(see eg \cite[Prop 2]{Pie} or \cite[(2.3)]{DP}).

The \textbf{greedy classes}
$\mathscr{G}_{q}^{\alpha}(\mathcal{B},\mathbb{B})$ are defined as
before replacing the role of $\sigma_N(x)$ by the error of greedy
approximation $\ga_N(x)$ given in \eqref{defnerrorAppGreedy}, that
is
\begin{eqnarray}\label{defnApprxSpacesGre}
\mathscr{G}_{q}^{\alpha}(\mathcal{B},\mathbb{B}) =
\Big\{x\in\mathbb{B}:\|x\|_{\mathscr{G}^{\alpha}_{q}}
\equiv\|x\|_{\mathbb{B}} +
\Big[\sum_{N=1}^{\infty}\big(N^{\alpha}\gamma_{N}(x;\mathcal{B},\mathbb{B})\big)^{q}
\frac{1}{N}\Big]^{\frac{1}{q}}<\infty\Big\}
\end{eqnarray}
(and similarly for $q = \infty$).
We also have the equivalence \Be
\|x\|_{\mathscr{G}_{q}^{\alpha}}\approx \|x\|_{\mathbb{B}} +
\Big[\sum_{k=0}^{\infty}\bigl(2^{k\alpha}\ga_{2^{k}}(x)\bigr)^{q}\Big]^{\frac{1}{q}},\label{Gaq-dyad}
\Ee since $\gamma_{N}(x)$ is non-increasing by the lattice
property in $\S$2.1.


Since $\sigma_{N}(x)\leq
\gamma_{N}(x)$ for all $x \in \mathbb{B}$
it is clear  that\footnote{Here, as in the rest of the paper,
$X\hookrightarrow Y$ means $X\subset Y$ and there exists $C>0$
such that $\|x\|_Y \leq C \|x\|_X$ for all $x\in X$. The equality
of spaces $X=Y$ is interpreted as $X\hookrightarrow Y$ and
$Y\hookrightarrow X$.}
\begin{eqnarray}\label{InclBesGreApprSpa}
\mathscr{G}_{q}^{\alpha}(\mathcal{B},\mathbb{B})\hookrightarrow
\mathcal{A}_{q}^{\alpha}(\mathcal{B},\mathbb{B})\,.
\end{eqnarray}
 When $\mathcal{B}$ is a greedy basis
in $\mathbb{B}$  it  holds that
$\mathscr{G}_{q}^{\alpha}(\mathcal{B},\mathbb{B}) =
\mathcal{A}_{q}^{\alpha}(\mathcal{B},\mathbb{B})$ with equivalent
quasi-norms. For non greedy bases, however, the inclusion may be
strict, and the classes $\mathscr{G}_{q}^{\alpha}$ may not even be
linear spaces (see section \ref{Noninclusions} below).

\subsection{Discrete Lorentz
Spaces}\label{subDiscLorSpac}

Let $\eta=\{\eta(k)\}_{k=1}^\infty$ be a sequence so that\Benu
\item[(a)] 
$0<\eta(k)\leq \eta(k+1)$ for all $k=1,2,\ldots$ and
$\lim_{k\longrightarrow\infty}\eta(k) = \infty$. \smallskip

\item[(b)] $\eta$ is \emph{doubling}, that is, $\eta(2k)\leq
C\eta(k)$ for all $k=1,2,\ldots,$ and some $C>0$. \Eenu We shall
denote the set of all such sequences by $\SW$. If $\eta\in\SW$ and
$0<r\leq\infty,$ the \textbf{weighted discrete Lorentz space}
{\boldmath$\ell_{\eta}^{r}$} is defined as
\[\ell_{\eta}^{r} = \Big\{\textbf{s}=\{s_k\}_{k=1}^\infty\in
\mathfrak{c}_{0}\mid
\|\textbf{s}\|_{\ell_{\eta}^{r}}\equiv\Big[\sum_{k=1}^{\infty}
\bigl(\eta(k)s^*_{k}\bigr)^{r}\tfrac{1}{k}\Big]^{\frac{1}{r}}<\infty\Big\}
\]
(with
$\|\textbf{s}\|_{\ell_{\eta}^{\infty}}=\sup_{k\in\SN}\eta(k)s^*_k$
when $r=\infty$). Here $\{s^*_k\}$ denotes the decreasing
rearrangement of $\{|s_k|\}$, that is $s^*_k=|s_{\pi(k)}|$ where
$\pi$ is any bijection of $\mathbb{N}$ such that $|s_{\pi(k)}|\geq
|s_{\pi(k+1)}|$ for all $k=1,2,\ldots$ (since we are assuming
 $\lim_{k\to\infty}s_k=0$ such $\pi$'s always exist).
 When $\eta\in\SW$ the set $\ell_{\eta}^{r}$ is a quasi-Banach
  space (see eg \cite[$\S$2.2]{CRS}).
Equivalent quasi-norms are given by
 \begin{equation} \label{dyadic}
  \|\textbf{s}\|_{\ell_\eta^r}\approx
\Big[\sum_{j=0}^\infty
 \big(\eta(\kappa^j) s^*_{\kappa^j}\big)^r\Big]^{1/r}\,,
 \end{equation}
for any fixed integer $\kappa>1$. Particular examples are the
classical Lorentz sequence spaces $\ell^{p,r}$ (with
$\eta(k)=k^{1/p}$), and the Lorentz-Zygmund spaces
$\ell^{p,r}(\log\ell)^\gamma$ (for which
$\eta(k)=k^{1/p}\log^\gamma(k+1)$; see eg \cite[p. 285]{BS}).

Occasionally we will need to assume a stronger condition on the
weights $\eta.$ For an increasing sequence $\eta$ we define
  $$M_{\eta}(m) =
\sup_{k\in\SN}\frac{\eta(k)}{\eta(mk)}, \quad m=1,2,3,\ldots.
  $$
Observe that we always have $M_{\eta}(m)\leq 1.$ We shall say that
$\eta\in\SW_+$ when  $\eta\in\SW$ and there exists some integer
$\kappa>1$ for which $M_{\eta}(\kappa)<1.$ This is equivalent to
say that the ``lower dilation index'' $i_\eta>0$, where we let
\[
i_\eta\equiv\sup_{m\geq1}\frac{\log M_\eta(m)}{-\log m}\,.
\]
For example, $\eta=\{k^\al\log^\beta(k+1)\}$ has $i_\eta=\alpha$,
and hence $\eta\in\SW_+$ iff $\al>0$. In general, if $\eta$ is
obtained from a increasing function $\phi:\SR^+\to\SR^+$ as
$\eta(k)=\phi(ak)$, for some fixed $a>0$, then $i_\eta>0$ iff
$i_\phi>0$, the latter denoting the standard lower dilation index
of $\phi$ (see eg \cite[p. 54]{KPS} for the definition).

Below we will need the following result:

\begin{lemma}\label{lemaLorSpa}
If $\eta\in\SW_+$ then there exists a constant $C> 0$ such that
\Be\sum_{j=0}^{n}\eta(\kappa^{j})\leq C\eta(\kappa^{n}),\quad
\forall\; n\in\SN,\label{SW+}\Ee
 where $\kappa>1$  is an integer as in the definition of
 $\SW_+.$
 \end{lemma}

 \begin{proof}  Write $\dt =
 M_{\eta}(\kappa)<1.$ By definition $M_{\eta}(\kappa) \geq
 {\eta(\kappa^{j})}/{\eta(\kappa^{j+1})}$, and therefore
 \begin{eqnarray}\label{proofLem}
 \eta(\kappa^{j})\leq \dt
 \eta(\kappa^{j+1}),\quad \forall\; j=0,1,2,\ldots.
 \end{eqnarray}
 Iterating (\ref{proofLem}) we deduce that $\eta(\kappa^{j})\leq
 \dt^{n-j}\eta(\kappa^{n}),$ for $j=0,1,2,\ldots,n$ and hence
 $$
 \sum_{j=0}^{n}\eta(\kappa^{j})\leq\eta(\kappa^{n})\sum_{j=0}^{n}\dt^{n-j}\leq
 \eta(\kappa^{n})\frac{1}{1-\dt}\,.
 $$
\end{proof}

\begin{Remark} \label{remark2.3}
{\rm If $\eta$ is increasing and doubling, then $
\{k^\alpha\,\eta(k)\}\in\SW_+$ for all $\al>0$. Also, if $\eta
\in\SW_+$ then $\eta^r \in\SW_+$, for all $r>0$.}
\end{Remark}

We now estimate  the \emph{fundamental function} of $\ell^r_\eta$.
We shall denote the indicator sequence of $\Ga\subset\SN$ by
$1_\Ga$, that is the sequence with entries 1 for $j\in\Ga$ and 0
otherwise.

\begin{lemma}\label{lemanorma}
 (a) If  $\eta\in\SW$ then
  \[\big\|1_\Ga\big\|_{\ell_\eta^\infty} =
  \eta(|\Ga|),\quad \forall\;\mbox{\rm finite }\Ga\subset\SN.\]

(b) If $\eta\in\SW_+$ and $r\in(0,\infty)$ then
\[\big\|1_\Ga\big\|_{\ell_\eta^r} \approx
  \eta(|\Ga|),\quad \forall\;\mbox{\rm finite }\Ga\subset\SN\]
  with the constants involved independent of $\Ga$.
 \end{lemma}

 \begin{proof}
Part (a) is trivial since $\eta$ is increasing. To prove (b) use
(\ref{dyadic}) and the previous lemma.
 \end{proof}

Finally, as mentioned in $\S1$, given a (normalized) basis $\cB$ in $\SB$ we shall consider the following
subspaces
\[
\ell^q_\eta(\cB,\SB):=\Bigl\{x=\sum_{j=1}^{\infty} c_j e_j \in\SB\;\colon\; \{c_j\}_{j=1}^\infty\in\ell^q_\eta\Bigr\},
\]
endowed with the quasi-norm
$\|x\|_{\ell^q_\eta(\cB,\SB)}:=\|\{c_j\}\|_{\ell^q_\eta}$. These spaces are not necessarily complete,
but they are when
\[
\|\sum_j c_j e_j\|_{\mathbb{B}}\leq C \|\{c_j\}\|_{\ell_{\eta}^q}, \quad \forall\mbox{ finite }\{c_j\},
\]
a property which holds in certain situations (see eg Remark
\ref{RemProofThJacIne1}). When this is the case, the space
$\ell^q_\eta(\cB,\SB)$ is just an isomorphic copy of $\ell^q_\eta$
inside $\SB$.

\subsection{Democracy Functions}\label{sbDemocFunc}

Following \cite{KT}, a (normalized) basis $\mathcal{B}$ in a
quasi-Banach space $\mathbb{B}$ is said to be \textbf{democratic}
if there exists $C>0$ such that
\[\Big\|\sum_{k\in\Gamma} e_{k}\Big\|_{\mathbb{B}}\leq
C\Big\|\sum_{k\in\Gamma'} e_{k}\Big\|_{\mathbb{B}},
\]
for all finite sets $\Gamma,\Gamma' \subset \mathbb{N}$ with the
same cardinality. This notion allows to characterize greedy bases
as those which are both unconditional and democratic \cite{KT}.

As we recall in $\S5$, wavelet bases are well known examples of
greedy bases for many function spaces, such as $L^p$, Sobolev, or
more generally, the Triebel-Lizorkin spaces.
However, they are not democratic in some other instances such as
$BMO$, or the  Orlicz $L^\Phi$ and Lorentz $L^{p,q}$ spaces (when
these are different from $L^p$). In fact, it is proved in
\cite{Wo1} that the Haar basis is democratic in a rearrangement
invariant space $\mathbb{X}$ in $[0,1]$ if and only if $\mathbb{X}
= L^{p}$ for some $p\in(1,\infty).$

Thus, non-democratic bases are also common. To quantify the
democracy of a (normalized) system $\cB=\{e_j\}_{j=1}^\infty$ in
$\SB$
one introduces the following concepts:
\[h_{r}(N;\mathcal{B},\mathbb{B})\equiv
\sup_{|\Gamma|=N}\Big\|\sum_{k\in\Gamma}e_{k}\Big\|_{\mathbb{B}}
\quad \mbox{and} \quad
 h_\ell(N;\mathcal{B},\mathbb{B})\equiv\inf_{|\Gamma|=
N}\Big\|\sum_{k\in\Gamma}e_{k}\Big\|_{\mathbb{B}}\,,
\]
which we shall call the {\bf right and left democracy functions
of} $\mathcal{B}$ (see also \cite{DKKT,KT2,GHM}).
We shall omit $\mathcal{B}$ or $\mathbb{B}$ when these are
understood from the context.

Some general properties of $h_\ell$ and $h_r$ are proved in the
next proposition.

\begin{proposition}\label{propPropertyRLDemoc}
Let $\cB=\{e_j\}_{j=1}^\infty$ be a (normalized) unconditional
basis in $\SB$ with the lattice property from $\S$2.1. Then \Benu
\item[(a)]$1\leq h_\ell(N)\leq h_r(N)\leq N^{1/\rho}$,
$\forall\;N=1,2,\ldots$, where $\rho=\rho_\SB$ is as in
\eqref{rhotrianIneq}. \item[(b)] $h_{\ell}(N)$ and $h_r(N)$ are
non-decreasing in $N=1,2,3\ldots$
 \item[(c)] $h_{r}(N)$ is doubling, that is, $\exists\; c>0$ such that $h_{r}(2N)\leq c\,
 h_{r}(N)$,
 $\forall\;N\in\SN$.
 \item[(d)]
  There exists $c\geq 1$
 such that $h_\ell(N+1)\leq c \,h_\ell(N)$ for all $N=1,2,3\ldots$
\Eenu \end{proposition}

 \begin{proof} (a) and (b) follow immediately from the lattice property of $\mathcal{B}$
 and the $\rho$-triangular inequality.\\
 (c) Given $N\in\mathbb{N},$ choose $\Gamma\subset\mathbb{N}$
 with $|\Gamma| = 2N$ such that
   $
   \big\|\sum_{k\in\Gamma}e_{k}\big\|_{\mathbb{B}}\geq h_{r}(2N)/2$.
 Partitioning arbitrarily $\Gamma = \Gamma'\cup\Gamma''$ with $|\Gamma'| =
 |\Gamma''|=N$,
 and using the $\rho$-power triangle inequality, one easily
 obtains
 \[
\tfrac{1}{2}h_{r}(2N)\leq\Big\|\sum_{k\in\Gamma}e_{k}\Big\|_{\mathbb{B}}
 = \Big\| \sum_{k\in\Gamma'}e_{k}+
 \sum_{k\in\Gamma''}e_{k}\Big\|_{\mathbb{B}}\leq 2^{1/\rho}h_{r}(N)\,.
 \]
(d) Given $N\in\mathbb{N},$  choose $\Gamma\subset\mathbb{N}$ with
$|\Gamma| = N$ such that
   $
   \big\|\sum_{k\in\Gamma}e_{k}\big\|_{\mathbb{B}}\leq 2h_\ell(N).
   $
Let $\Gamma' = \Gamma\cup\{k_{o}\}$ for any  $k_{o}\notin \Gamma$.
Then
 \begin{eqnarray*}
 h_\ell(N+1)\leq \Big\|\sum_{k\in\Gamma'}e_{k}\Big\|_{\mathbb{B}}\leq
 \Big(\Big\|\sum_{k\in\Gamma}e_{k}\Big\|_{\mathbb{B}}^{\rho} + 1\Big)^{1/\rho}
 \leq (2^{\rho}[h_\ell(N)]^{\rho}+1)^{1/\rho}.
 \end{eqnarray*}
 Thus, using (a) we obtain
   $   h_\ell(N+1)\leq (2^{\rho}+1)^{\frac{1}{\rho}}h_\ell(N)\leq 2\cdot2^{1/\rho}h_\ell(N).$
 \end{proof}


\begin{Remark}
{\rm We do not know whether property (d) can be improved to show
that $h_\ell(N)$ is actually doubling. This seems however to be
case in all the examples we have considered below (see $\S5$).}
\end{Remark}

\section{Right Democracy and Jackson Type
Inequalities}\label{secRigtDemJackIneq}

Our first result deals with inclusions for the greedy classes
$\GaqB$.

\begin{theorem}\label{ThGreJacksonIneq}
Let $\mathcal B =\{e_{j}\}_{j=1}^\infty$ be a (normalized)
unconditional basis in $\SB$. Fix  $\alpha > 0$ and $q\in(0,
\infty)$. Then, for any sequence $\eta$ such that
$\{k^\al\eta(k)\}_{k=1}^\infty\in\SW_+$ the following statements
are equivalent:

\sline 1. There exists $C> 0$ such that for all $N=1,2,3,\ldots$
    \begin{eqnarray}\label{ThJacksonIneq1}
 \Big\|\sum_{k\in\Gamma}e_{k}\Big\|_{\mathbb{B}}\leq
 C\eta(N)\,,\quad\forall\;\Ga\subset\SN\mbox{ with }|\Ga|=N.
     \end{eqnarray}

\sline 2. Jackson type inequality for
 $\ell_{k^{\alpha}\eta(k)}^{\,\infty}(\mathcal{B},\mathbb{B})$:
 $\exists\;C_{\alpha}>0$ such
 that $\forall\;N=0,1,2\ldots$
 \begin{eqnarray}\label{ThJacksonIneq2}
 \gamma_{N}(x)\leq
 C_{\alpha}(N+1)^{-\alpha}\|x\|_{\ell_{k^{\alpha}\eta(k)}^{\,\infty}(\mathcal{B},\mathbb{B})},
 \quad \forall\; x \in
 \ell_{k^{\alpha}\eta(k)}^{\,\infty}(\mathcal{B},\mathbb{B}).
 \end{eqnarray}

\sline 3.  $
     \ell_{k^{\alpha}\eta(k)}^{\infty}(\mathcal{B},\mathbb{B})\hookrightarrow
 \mathscr{G}_{\infty}^{\alpha}(\mathcal{B},\mathbb{B})\,.
     $

\sline 4.
$\ell_{k^{\alpha}\eta(k)}^{q}(\mathcal{B},\mathbb{B})\hookrightarrow
 \mathscr{G}_{q}^{\alpha}(\mathcal{B},\mathbb{B})\,.$

\sline 5. Jackson type inequality for
$\ell_{k^{\alpha}\eta(k)}^{\,q}(\mathcal{B},\mathbb{B})$:
   $\exists\;C_{\alpha,q}>0$ such that $\forall\;N
=0,1,2,\ldots$
 \begin{eqnarray}\label{ThJacksonIneq3}
 \gamma_{N}(x)\leq
C_{\alpha,q}(N+1)^{-\alpha}\|x\|_{\ell_{k^{\alpha}\eta(k)}^{\,q}(\mathcal{B},\mathbb{B})},
\quad \forall\;x \in
 \ell_{k^{\alpha}\eta(k)}^{\,q}(\mathcal{B},\mathbb{B})\,.
\end{eqnarray}

\end{theorem}

 \begin{proof}
 ``$1\Rightarrow 2$''
 Let $x= \sum_{k\in\mathbb{N}}c_{k}e_{k}\in \ell_{k^{\alpha}\eta(k)}^{\infty}
 (\mathcal{B},\mathbb{B})$
and let $\pi$ be a bijection of $\mathbb{N}$ such that
\begin{eqnarray}\label{nonincRearraThJ}
|c_{\pi(k)}| \geq |c_{\pi(k+1)}|, \quad k=1,2,3,\ldots
\end{eqnarray}
For fixed $N=0,1,2,\ldots$, denote $\la_j=2^j(N+1)$. Then, the
$\rho$-power triangle inequality and (\ref{LaticeIneq}) give
\begin{eqnarray*}
\bigl\|x -G_{N}^{\pi}(x)\bigr\|_{\mathbb{B}}^{\rho} &=&
\Big\|\sum_{k=N+1}^{\infty}c_{\pi(k)}e_{\pi(k)}\Big\|_{\mathbb{B}}^{\rho}
 \,\leq\, \sum_{j=0}^{\infty}\Big\|{\Ts\sum_{\la_j\leq k
<\la_{j+1}}c_{\pi(k)}e_{\pi(k)}}\Big\|_{\mathbb{B}}^{\rho}\,\\
&\leq &\sum_{j=0}^{\infty} |c_{\pi(\la_j)}|^{\rho}\,
\Big\|{\Ts\sum_{\la_j\leq k <
\la_{j+1}}e_{\pi(k)}}\Big\|_{\mathbb{B}}^{\rho}.
\end{eqnarray*}
There are exactly $\la_j=2^{j}(N+1)$ elements in the interior sum,
so using \eqref{ThJacksonIneq1} we obtain \Bea \|x-G_{N}^{\pi}(x)\|_{\mathbb{B}}^{\rho}
&\leq& C^{\rho}\sum_{j=0}^{\infty}\bigl(c^*_{\la_j}
\eta(\la_j)\bigr)^{\rho}\,=\,
 C^{\rho}\sum_{j=0}^{\infty}\bigl(\la_j^{\alpha}c^*_{\la_j}\eta(\la_j)
 \bigr)^{\rho}\,\la_j^{-\alpha\rho}
 \nonumber\\ &\leq&
 C^{\rho}\|x\|_{\ell_{k^{\alpha}\eta(k)}^{\infty}(\mathcal{B},\mathbb{B})}^{\rho}
 \,\mbox{\small$(N+1)^{-\alpha\rho}\,
 \sum_{j=0}^{\infty}{2^{-j\alpha\rho}}$} \nonumber\\ &=&
 C_{\alpha, \rho}\,(N+1)^{-\alpha\rho}\,\|x\|_{\ell_{k^{\alpha}\eta(k)}^{\infty}(\mathcal{B},\mathbb{B})}^{\rho}
 \,. \nonumber
 \Eea
 The result follows taking the supremum over all bijections $\pi$
 satisfying (\ref{nonincRearraThJ}).

 \begin{Remark}\label{RemProofThJacIne1} {\rm The  special case $N=0$
 in \eqref{ThJacksonIneq2} says that
 \Be\|x\|_{\mathbb{B}}\leq C \|x\|_{\ell_{k^{\alpha}\eta(k)}^{\infty}
 (\mathcal{B},\mathbb{B})},\label{aux1}\Ee
which in particular implies
 $\ell_{k^{\alpha}\eta(k)}^{q}(\mathcal{B},\mathbb{B})\hookrightarrow
 \mathbb{B}$, for all $q\in(0,\infty]$.
}
 \end{Remark}

\noindent ``$2\Rightarrow 3$'' This is immediate from the definition of
 $\mathscr{G}_{\infty}^{\alpha}$ (and Remark \ref{RemProofThJacIne1}), since
   $$
  \|x\|_{\mathscr{G}_{\infty}^{\alpha}(\mathcal{B},\mathbb{B})}:=\|x\|_\SB+ \sup_{N\geq
 1}N^{\alpha}\gamma_{N}(x)\leq
 C_{\alpha}\|x\|_{\ell_{k^{\alpha}\eta(k)}^{\infty}(\mathcal{B},\mathbb{B})}.
   $$

\noindent ``$3\Rightarrow 1$'' Let
 $\Gamma \subset\mathbb{N}$ with $|\Gamma| = N.$ Choose $\Gamma'$
 with $|\Gamma'| = N$ and so that $\Gamma\cap \Gamma' =
 \emptyset$, and  consider
   $
   x =
 \sum_{k\in \Gamma} e_{k} +
 \sum_{k\in \Gamma'} 2e_{k}\,.
 $ Then
 \begin{eqnarray}\label{proofThJacIneq1}
 \gamma_{N}(x)=
 \big\|\sum_{k\in\Gamma}e_{k}\big\|_{\mathbb{B}}\,,
 \end{eqnarray}
 and therefore
\Be
 N^{\alpha}\big\|\sum_{k\in\Gamma}e_{k}\big\|_{\mathbb{B}}\,=\, N^{\alpha}\gamma_{N}(x) \,\leq\,
 \|x\|_{\mathscr{G}_{\infty}^{\alpha}(\mathcal{B},\mathbb{B})}.
\label{J2a}\Ee
 On the other hand, call $\om(k)=k^\al \eta(k)$.
 By monotonicity, Lemma \ref{lemanorma} and the doubling property of $\om$ we have
 \Be
 \|x\|_{\ell_{\om}^{\infty}(\mathcal{B},\mathbb{B})}
 \leq 2\big\|1_{\Gamma \cup \Gamma'}\big\|_{\ell_{\om}^{\infty}}
 = 2 \om(2N) \leq\, c \,\om(N)\,.
 \label{J2b}\Ee
 Combining \eqref{J2a} and \eqref{J2b} with the inclusion $\ell_{k^{\alpha}\eta(k)}^{\infty}(\mathcal{B},\mathbb{B})\hookrightarrow
 \mathscr{G}_{\infty}^{\alpha}(\mathcal{B},\mathbb{B})$ gives \eqref{ThJacksonIneq1}.

\bline ``$5\Rightarrow 1$'' Let $\Gamma \subset \mathbb{N}$ with
$|\Gamma| =
 N$, and choose $\Gamma'$ and $x$ as in the proof of $3\Rightarrow
 1.$ As before call $\om(k)=k^\al \eta(k)$.
Then Lemma \ref{lemanorma} and the assumption $\om\in\SW_+$ give
 $$
 \|x\|_{\ell_{\om}^{q}(\mathcal{B},\mathbb{B})}
 \leq 2\big\|1_{\Gamma \cup \Gamma'}\big\|_{\ell_{\om}^{q}}
 \approx  \om(2N) \leq\, c \,\om(N)\,.
 $$
Since we are assuming 5 we can write (recall
(\ref{proofThJacIneq1}))
 $$
 \big\|\sum_{k\in\Gamma}{e_{k}}\big\|_{\mathbb{B}}
 = \gamma_{N}(x)
 \leq
 C_{\alpha,\rho}(N+1)^{-\alpha}\|x\|_{\ell_{\om}^{q}(\mathcal{B},\mathbb{B})}
 \lesssim N^{-\al}\om(N)=\eta(N),
 $$
which proves \eqref{ThJacksonIneq1}.

\bline ``$1\Rightarrow 4$''
The proof is similar to $1\Rightarrow2$ with a few modifications we indicate next.
Given $x\in
 \ell_{k^{\alpha}\eta(k)}^{q}(\mathcal{B},\mathbb{B})$ and $\pi$ as in \eqref{nonincRearraThJ}
we write $x=\sum_{j=-1}^\infty\sum_{2^j<
k\leq2^{j+1}}c_{\pi(k)}e_{\pi(k)}$. Then
arguing as before (with $N=2^m$) we obtain
\[ \|x -  G_{2^m}^{\pi}(x)\|_{\mathbb{B}}^{\mu}\,\leq
 \,\sum_{j=m}^{\infty} |c_{\pi(2^j)}|^{\mu}\,
\Big\|{\Ts\sum_{2^j< k \leq
2^{j+1}}e_{\pi(k)}}\Big\|_{\mathbb{B}}^{\mu},
\]
where
we choose now any $\mu<\min\{q,\rho_\SB\}$.  Taking the supremum over all
 $\pi$'s and using \eqref{ThJacksonIneq1} we obtain
\[\gamma_{2^m}(x;\mathcal{B},\mathbb{B})^{\mu} \leq
C^{\mu}\sum_{j=m}^{\infty}
\bigl(c^*_{2^j}\,\eta(2^j)\bigr)^{\mu}.
\]
Therefore
$$
  \Big[\sum_{m=0}^{\infty}\bigl(2^{m\alpha}\gamma_{2^m}(x)\bigr)^{q}\Big]^{\frac{1}{q}}  \leq
C\Big[\sum_{m=0}^{\infty}2^{m\alpha
q}\Bigl(\sum_{j=0}^{\infty}\bigl[c^*_{2^{j+m}}\,\eta(2^{j+m})\bigr]^{\mu}\Bigr)^{q/\mu}\Big]^{1/q}.
$$
Since $q/\mu > 1,$ we can use Minkowski's inequality on the right hand side to obtain
\Beas
\Big[\sum_{m=0}^{\infty}\bigl(2^{m\alpha}\gamma_{2^m}(x)\bigr)^{q}\Big]^{\frac{1}{q}}  &\leq &
C\Big[\sum_{j=0}^{\infty}\Bigl(\sum_{m=0}^{\infty}2^{m\alpha
q}\bigl[c^*_{2^{j+m}}\,\eta(2^{j+m})\bigr]^q\Bigr)^{\mu/q}\Big]^{1/\mu}
\\&= &
C\Big[\sum_{j=0}^{\infty}2^{-j\al\mu}\Bigl(\sum_{\ell=j}^{\infty}2^{\ell\alpha
q}\bigl[c^*_{2^\ell}\,\eta(2^\ell)\bigr]^q\Bigr)^{\mu/q}\Big]^{1/\mu}
\, \leq \, C'\,\|\{c_k\}\|_{\ell_{k^{\alpha}\eta(k)}^{\,q}}.
\Eeas
This implies the desired estimate
\[
\|x\|_{\mathscr{G}_{q}^{\alpha}(\mathcal{B},\mathbb{B})}\,\lesssim\,\|\{c_k\}\|_{\ell_{k^{\alpha}\eta(k)}^{\,q}},\]
using the dyadic expressions for the norms in \eqref{Gaq-dyad} and \eqref{dyadic}
(and  Remark \ref{RemProofThJacIne1}).

\bline ``$4\Rightarrow 5$'' This is trivial since 4 implies
 $\ell_{k^{\alpha}\eta{k}}^{q}(\mathcal{B},\mathbb{B})
\hookrightarrow\mathscr{G}_{q}^{\alpha}(\mathcal{B},\mathbb{B})
\hookrightarrow\mathscr{G}_{\infty}^{\alpha}(\mathcal{B},\mathbb{B})$,
and this clearly gives \eqref{ThJacksonIneq3}.
\end{proof}

\begin{Remark}
{\rm The equivalences 1 to 3 remain true under the weaker
assumption $\{k^\al\eta(k)\}\in\SW$.}
\end{Remark}

\begin{Remark} \label{RemProofThJacIne2}
{\rm Observe that if any of the statements in 2 to 5 of Theorem
\ref{ThGreJacksonIneq} holds for one fixed $\alpha >0$ and
$q\in(0,\infty]$, then the assertions remain true for all $\alpha$
and $q$ (as long as $\{k^\al\eta(k)\}\in\SW_+$), since the
statement in 1 is independent of these parameters.}
\end{Remark}

\begin{corollary}
\label{CorProofThJacIne2} {\bf{Optimal inclusions into
\boldmath{$\mathscr{G}_{q}^{\alpha}$}.}}

\noindent Let $\mathcal B$ be a (normalized) unconditional basis
in $\SB$. Fix  $\al>0$ and $q\in(0,\infty]$. Then \Be
\ell_{k^{\alpha}h_r(k)}^{q}(\mathcal{B},\mathbb{B})
\hookrightarrow\mathscr{G}_{q}^{\alpha}(\mathcal{B},\mathbb{B}).
\label{low_inc_G}\Ee Moreover, if $\om\in\SW_+$ then,
$\ell_{\om}^{q}(\mathcal{B},\mathbb{B})
\hookrightarrow\mathscr{G}_{q}^{\alpha}(\mathcal{B},\mathbb{B})$
if and only if $\om(k)\gtrsim k^\al h_r(k)$.
\end{corollary}
\begin{proof}
For $ q<\infty$, the inclusion \eqref{low_inc_G} is an application
of 4 in the theorem with $\eta=h_r$ (after noticing that $\{k^\al
h_r(k)\}\in\SW_+$ by Proposition \ref{propPropertyRLDemoc} and
Remark \ref{remark2.3}). The second assertion  is just a restatement
of $1\Leftrightarrow4$ with $\eta(k)=\om(k)/k^\al$. For $q=\infty$
use 3 instead of 4.
 \end{proof}

We now prove similar results for the approximation spaces
$\AaqB$.

\begin{theorem}\label{ThBestJackIneq}
Let $\mathcal B =\{e_{j}\}_{j=1}^\infty$ be a (normalized)
unconditional basis in $\SB$. Fix  $\alpha > 0$ and $q\in(0,
\infty]$. Then, for any sequence $\eta\in\SW_+$ the following are
equivalent:

\sline 1. There exists $C> 0$ such that for all $N=1,2,3,\ldots$
\begin{eqnarray}\label{ThBestJackIneq1}
\Big\|\sum_{k\in\Gamma}{e_{k}}\Big\|_{\mathbb{B}}\leq
C\eta(N)\,,\quad\forall\;\Ga\subset\SN\mbox{ with }|\Ga|=N.
\end{eqnarray}

%

\sline 2.
$\ell_{k^{\alpha}\eta(k)}^{q}(\mathcal{B},\mathbb{B})\hookrightarrow
\mathcal{A}_{q}^{\alpha}(\mathcal{B},\mathbb{B})$ .

\sline 3. Jackson type inequality for
$\ell_{k^{\alpha}\eta(k)}^{\,q}(\mathcal{B},\mathbb{B})$:
$\exists\;C_{\alpha, q} >0$ such that $\forall\;N= 0,1, 2, \dots$
\begin{eqnarray}\label{ThBestJackIneq3}
\sigma_{N}(x)\leq C_{\alpha,
q}(N+1)^{-\alpha}\|x\|_{\ell_{k^{\alpha}\eta(k)}^{\,q}(\mathcal{B},\mathbb{B})},
\quad \forall\;x \in
 \ell_{k^{\alpha}\eta(k)}^{q}(\mathcal{B},\mathbb{B})\,.
\end{eqnarray}
\end{theorem}

\begin{proof} $1\Rightarrow 2$ follows directly from Theorem
\ref{ThGreJacksonIneq} and $\Gaq\hookrightarrow\mathcal{A}^\al_q$.
 Also, $2\Rightarrow 3$ is trivial since
$\mathcal{A}_q^{\alpha}\hookrightarrow\mathcal{A}_{\infty}^{\alpha}$,
and 3 is equivalent to
{\small$\ell_{k^{\alpha}\eta(k)}^{q}(\mathcal{B},\mathbb{B})\hookrightarrow
\mathcal{A}_{\infty}^{\alpha}$}.

\sline We must show $3\Rightarrow1$. Let  $\kappa>1$ be a fixed
integer
 as in the definition of the class $\SW_+$ (and in particular
 satisfying  \eqref{SW+}), and denote $1_{\Dt} = \sum_{k\in\Dt}{e_{k}}$ for
 a set $\Dt\subset\SN$. For any $\Ga_n\subset\SN$ with $|\Ga_n|=\kappa^n$,
 we can find a subset $\Gamma_{n-1}$ with $|\Gamma_{n-1}| = \kappa^{n-1}$ such that
  $$
  \|1_{\Gamma_{n}}- 1_{\Gamma_{n-1}}\|_{\mathbb{B}}\leq
2\sigma_{\kappa^{n-1}}(1_{\Gamma_{n}}).
  $$
Repeating this argument we choose $\Gamma_{j-1}\subset\Ga_{j}$ with $|\Gamma_{j}| = \kappa^{j}$ and so that
$$
  \|1_{\Gamma_{j}}- 1_{\Gamma_{j-1}}\|_{\mathbb{B}}\leq
2\sigma_{\kappa^{j-1}}(1_{\Gamma_{j}}), \quad\mbox{for }
j=1,2\ldots,n\,.
  $$
Setting $\Gamma_{-1} =
\emptyset,$ and using the $\rho$-power triangle inequality we see that
$$
\|1_{\Gamma_{n}}\|_{\mathbb{B}}^{\rho} = \Big\|\sum_{j=0}^{n}
1_{\Gamma_{j}} - 1_{\Gamma_{j-1}}\Big\|_{\mathbb{B}}^{\rho}
\leq\sum_{j=0}^{n}\| 1_{\Gamma_{j}}-
1_{\Gamma_{j-1}}\|_{\mathbb{B}}^{\rho} \leq 2^\rho \sum_{j=0}^n
\sigma_{\kappa^{j-1}}(1_{\Gamma_{j}})^\rho\,.
$$
Now, the hypothesis \eqref{ThBestJackIneq3} and Lemma \ref{lemanorma} give
 \[\sigma_{\kappa^{j-1}} (1_{\Gamma_{j}})\, \lesssim\,
 \kappa^{-j\alpha}
\|1_{\Gamma_{j}}\|_{\ell_{k^{\alpha}\eta(k)}^{\,q}(\mathcal{B},\mathbb{B})}\,\approx\,
\eta(\kappa^j).\]
 Thus, combining these two expressions we obtain
\Be
\|1_{\Gamma_{n}}\|_{\mathbb{B}} \,\lesssim\,\Bigl[ \sum_{j=0}^n \eta(\kappa^j)^\rho \Bigr]^{1\rho}
\,\leq C\, \eta(\kappa^n)\,, \label{last}\Ee
where the last inequality follows from the assumption $\eta\in\SW_+$ and Lemma \ref{lemaLorSpa}.
This shows \eqref{ThBestJackIneq1} when $N=\kappa^n$, $n=1,2,\ldots$
The general case follows easily using the doubling property of $\eta$.
\end{proof}

\begin{Remark} \label{RemProofThJacIne2}
{\rm As before, if any of the statements in 2 or 3 holds for one
fixed $\alpha >0$ and $q\in(0,\infty]$, then the assertions remain
true for all $\alpha$ and $q$, since  1 is independent of these
parameters.}
\end{Remark}

\begin{Remark}{\rm
Observe also that $1\Rightarrow2\Rightarrow3$ hold with the weaker
assumption $\{k^\al\eta(k)\}\in\SW_+$ from Theorem
\ref{ThGreJacksonIneq} (and in particular hold for $\eta=h_r$ as
stated in \eqref{inclusions-general}). However, the stronger
assumption $\eta\in\SW_+$ is crucial to obtain $3\Rightarrow1$,
and cannot be removed as shown in Example 5.6 below. }\end{Remark}

\begin{corollary}
\label{CorProofThJacIne2} {\bf{Optimality of the inclusions into
\boldmath{$\cA_{q}^{\alpha}$}.}}

\noindent Let $\mathcal B$ be a (normalized) unconditional basis
in $\SB$. Fix  $\al>0$ and $q\in(0,\infty]$. Then \Be
\ell_{k^{\alpha}h_r(k)}^{q}(\mathcal{B},\mathbb{B})
\hookrightarrow\cA_{q}^{\alpha}(\mathcal{B},\mathbb{B}).
\label{low_inc_A}\Ee If for some $\om\in\SW_+$ we have
$\ell_{\om}^{q}(\mathcal{B},\mathbb{B})
\hookrightarrow\cA_{q}^{\alpha}(\mathcal{B},\mathbb{B})$, then
necessarily $\om(k)\gtrsim k^\al$. Moreover if
$\om(k)=k^\al\eta(k)$, with $\eta$ increasing and doubling, then

\sline (a) if $i_\eta>0$, then necessarily $\eta(k)\gtrsim
h_r(k)$, and hence
$\ell^q_{k^\al\eta(k)}\hookrightarrow\ell^q_{k^{\alpha}h_r(k)}$.

\sline (b) if $i_\eta=0$, then
 $\eta(k)\gtrsim h_r(k)/(\log k)^{1/\rho}$ and
$\ell^q_{k^\al\eta(k)}\hookrightarrow\ell^q_{\{k^{\alpha}h_r(k)/(\log k)^{1/\rho}\}}$.

\end{corollary}
\begin{proof}
The inclusion \eqref{low_inc_A} is actually a consequence of
\eqref{low_inc_G}. Assertion (a) is just
$2\Rightarrow3\Rightarrow1$ in the theorem. For assertion (b)
notice that in the last step of the proof of $3\Rightarrow1$, the
right hand inequality of \eqref{last} can always be replaced by
\[
\|1_{\Gamma_{n}}\|_{\mathbb{B}} \,\lesssim\,\Bigl[ \sum_{j=0}^n
\eta(\kappa^j)^\rho \Bigr]^{1\rho} \,\lesssim\,
\eta(\kappa^n)\,n^{1/\rho}
\]
when $\eta$ is increasing. Thus $h_r(N)\lesssim \eta(N)(\log
N)^{1/\rho}$ holds for $N=\kappa^n$, and by the doubling property
also for all $N\in\SN$. Finally, if
$\ell_{\om}^{q}(\mathcal{B},\mathbb{B})
\hookrightarrow\cA_{q}^{\alpha}(\mathcal{B},\mathbb{B})$ for some
general $\om\in\SW_+$, then given $\Ga\subset\SN$ with $|\Ga|=N$ we
trivially have\[
\om(N)\approx\|1_\Ga\|_{\ell^q_\om}\gtrsim\|1_\Ga\|_{\cA^\al_\infty}\geq
(N/2)^\al\,\sigma_{N/2}(1_\Ga)\,\geq \,(N/2)^\al.
\]

 \end{proof}
\BR Assertion (b) shows that the inclusion in \eqref{low_inc_A} is
optimal, except perhaps for a logarithmic loss. The logarithmic loss
may actually happen, as there are Banach spaces $\SB$ with
$h_r(N)\approx \log N$ and so that
\[
\mathcal{A}^\al_q(\SB)\,=\,\ell^{\,q}_{k^\al}\,=\,\ell^q_{\{k^{\alpha}h_r(k)/\log
k\}}.
\]
See Example 5.6 below.
\ER

\section{Left Democracy and Bernstein Type
Inequalities}\label{secLeftDemBernsIneq}

It is well known that upper inclusions for the approximation
spaces $\cA^\al_q$, as in \eqref{inclusions-general}, depend upon
Bernstein type inequalities. In this section we show how the left
democracy function of $\cB$ is linked with these two properties.

We first remark that, for each $\alpha > 0$ and $0 < q \leq
\infty$, the approximation classes $\cA^\al_q$ and $\Gaq$
satisfy trivial Bernstein inequalities, namely, there exists $C_{\al,q}>
0$ such that
 \begin{equation} \label{4.1}
 \|x\|_{\mathcal A_q^\alpha(\mathcal B, \mathbb B)} \leq
 \|x\|_{\mathscr G_q^\alpha(\mathcal B, \mathbb B)} \leq \,C_{\al,q}\, N^\alpha
 \|x\|_{\mathbb B},\quad\forall\;x\in \Sigma_N,\;\; N=1,2, \dots
 \end{equation}
This follows easily from the definition of the norms and the trivial estimates $\sigma_N(x)\leq\ga_N(x)\leq\|x\|_\SB$.

 We  start with a preliminary result which is essentially known in the literature (see eg
 \cite{Pie}). As usual $\cB=\{e_j\}_{j=1}^\infty$ is a fixed
 (normalized) unconditional basis in $\SB$.

\begin{proposition}\label{propBernIneq}
Let $\SE$ be a subspace of $\SB$, endowed with a quasi-norm
$\|.\|_{\mathbb E}$ satisfying the $\rho$-triangle inequality for
some $\rho=\rho_\SE$. For each $\alpha
>0$ the following are equivalent:

\sline 1. $\exists\;C_{\alpha}
> 0$ such that
$ \|x\|_{\mathbb{E}}\leq
\,C_{\alpha}\,N^{\alpha}\,\|x\|_{\mathbb{B}}$,
$\forall\;x\in\Sigma_{N},\ N=1,2,\ldots\,$

\sline 2. $\mathcal A_\rho^\alpha (\mathcal B, \mathbb B)
\hookrightarrow \mathbb E\,.$

\sline 3. $\mathscr G_\rho^\alpha (\mathcal B, \mathbb B)
\hookrightarrow \mathbb E\,.$
\end{proposition}

\begin{proof} ``$1\Rightarrow 2$'' Given
$x\in\mathcal{A}_{\rho}^{\alpha}(\mathcal{B},\mathbb{B}),$ by the
representation theorem for approximation spaces \cite{Pie} one can
write $ x = \sum_{k=0}^{\infty}x_{k}$ with $x_{k}
\in\Sigma_{2^{k}},\ k=0,1,2,\ldots,$ such that
\[\Big(\sum_{k=0}^{\infty}2^{k\alpha\rho}
\|x_{k}\|_{\mathbb{B}}^{\rho}\Big)^{1/\rho}\leq
C\|x\|_{\mathcal{A}_{\rho}^{\alpha}(\mathcal B, \mathbb B)}\,.
\]
The hypothesis 1 and the $\rho_\SE$-triangular inequality then
give
$$
\|x\|_{\mathbb E}^{\rho} \leq
\sum_{k=0}^{\infty}\|x_{k}\|_{\mathbb E}^{\rho}
\leq\,C_{\alpha}^{\rho}\,
\sum_{k=0}^{\infty}2^{k\alpha\rho}\|x_{k}\|_{\mathbb{B}}^{\rho}\,\leq
\,C'\,\|x\|_{\mathcal{A_{\rho}^{\alpha}}(\mathcal B, \mathbb
B)}^{\rho}.
$$

\sline ``$2 \Rightarrow 3$''. This follows from the trivial
inclusion $\mathscr{G_{\rho}^{\alpha}}(\mathcal B, \mathbb B)
\hookrightarrow \mathcal{A_{\rho}^{\alpha}}(\mathcal B, \mathbb
B)\,.$

\sline ``$3\Rightarrow 1$''. This is immediate using (\ref{4.1}).
\end{proof}

\begin{theorem}\label{ThBernsIneq}
Let $\cB=\{e_j\}_{j=1}^\infty$ be a
 (normalized) unconditional basis in $\SB$.
Fix  $\alpha > 0$ and $q\in(0, \infty]$. Then, for any increasing
and doubling sequence $\{\eta(k)\}$
the following statements are equivalent:

\sline 1. There exists $C> 0$ such that for all $N=1,2,3,\ldots$
  \begin{eqnarray} \label{left}
 \Big\|\sum_{k\in\Gamma}e_{k}\Big\|_{\mathbb{B}}\,\geq \,\tfrac1C\,\eta(N),\quad\forall\;\Ga\subset\SN\mbox{ with }|\Ga|=N.
  \end{eqnarray}
\noindent 2. Bernstein type inequality for
$\ell_{k^{\alpha}\eta(k)}^{\,q}(\mathcal{B},\mathbb{B})$:
$\exists\;C_{\alpha,q} > 0$ such that
  \Be\|x\|_{\ell_{k^{\alpha}\eta(k)}^{\,q}(\mathcal{B},\mathbb{B})}\,\leq\,
C_{\alpha,q}\,N^\al\,\|x\|_{\mathbb{B}}, \quad \forall\;x
\in\Sigma_{N},\ N=1,2,3,\dots
\label{ber1}\Ee
 \noindent3.
  $\mathcal{A}_{q}^{\alpha}(\mathcal{B},\mathbb{B})\hookrightarrow
\ell_{k^{\alpha}\eta(k)}^{\,q}(\mathcal{B},\mathbb{B})\,.
  $
  \sline 4.
$\mathscr{G}_{q}^{\alpha}(\mathcal{B},\mathbb{B})\hookrightarrow
\ell_{k^{\alpha}\eta(k)}^{\,q}(\mathcal{B},\mathbb{B}).$
\end{theorem}

\begin{proof}
``$1\Rightarrow2$''. Let $x=
\sum_{k\in\Gamma}c_{k}e_{k}\in\Sigma_{N}$. For any bijection  $\pi$
with $|c_{\pi(k)}|$ decreasing, and any integer $m\in\{1,\ldots, N\}$ we have
 $$
|c_{\pi(m)}|\,\eta(m) \leq C\,
|c_{\pi(m)}|\,\big\|\sum_{j=1}^{m}e_{\pi(j)}\big\|_{\mathbb{B}}
\leq
C\,\big\|\sum_{j=1}^{m}c_{\pi(j)}e_{\pi(j)}\big\|_{\mathbb{B}}\leq
C\|x\|_{\mathbb{B}}\,,
 $$
using  (\ref{LaticeIneq}) in the second inequality. This gives
 $$
\|x\|_{\ell_{k^{\alpha}\eta(k)}^{q}} = \Big[\sum_{m=1}^{N}
(m^{\alpha}\eta(m)c^*_{m})^q\frac{1}{m}\Big]^{1/q} \leq
C\|x\|_{\mathbb{B}}\Big[\sum_{m=1}^{N}m^{\al q}\frac{1}{m}\Big]^{1/q}\,\approx\,
\|x\|_{\mathbb{B}}\,N^{\alpha}.
 $$

\sline ``$2\Rightarrow 1$''. For any $\Gamma \subset \mathbb N$ with $|\Ga|=N$,
applying \eqref{ber1} to  $1_\Gamma = \sum_{k\in \Gamma} e_k$ we obtain
 $$
  \|1_\Gamma\|_{\mathbb B} \geq \,\tfrac{1}{C_{\al,q}}\,
  N^{-\alpha}\|1_\Gamma\|_{\ell_{k^\alpha \eta(k)}^q (\mathcal B,
  \mathbb B)} \,\gtrsim\, \eta(N),
 $$
where in the last inequality we have used $\|1_\Gamma\|_{\ell_{\om}^{\,q}}\gtrsim \om(N)$, when $\om\in\SW$.

\sline ``$2\Rightarrow 3$''. We have already proved that $1
\Leftrightarrow 2$; since 1 does not depend on $\alpha,q$, then $2$
actually holds for all 
$\tilde{\alpha}>0$. In particular, from Proposition
\ref{propBernIneq}, we have
\Be\cA^{\tilde\al}_\rho\hookrightarrow\mathbb E :=
\ell_{k^{\tilde\alpha} \eta(k)}^q (\mathcal B,\mathbb
B)\label{AE}\Ee for $\tilde{\alpha} \in(\frac\al2,\frac{3\al}2)$
and some sufficiently small $\rho>0$. Now, from the general theory
developed in \cite{DP}, the spaces $\mathcal{A}_{q}^{\alpha}$
satisfy a reiteration theorem for the real interpolation method,
and in particular
\begin{equation}\label{proofTHBernIneq}
\mathcal{A}_{q}^{\alpha}\,=\,
\big(\mathcal{A}_{q_0}^{\alpha_{0}},\mathcal{A}_{q_1}^{\alpha_{1}}\big)_{{1}/{2},\,q}\;,
\end{equation}
when $\alpha =(\alpha_{0} + \alpha_{1})/{2}$ with $\alpha_{1}>\alpha_{0}>0$,
and $q_0,q_1,q\in(0,\infty]$.
On the other hand, for the family of weighted Lorentz spaces it is known that
\begin{eqnarray}\label{proofTHBernIneq1}
\big(\ell_{\om_0}^{q},\ell_{\om_1}^{q}\big)_{\theta,\,q}\,
= \,\ell_{\om}^{q}\,,\quad \mbox{ {\small$0<\theta<1,\quad 0<q\leq\infty,$}}
\end{eqnarray}
when $\om_0,\om_1\in\SW_+$ and $\om=\om_0^{1-\theta}\om_1^\theta$ (see eg \cite[Theorem 3]{MerApp}).
Thus, for fixed $\al$ and $q$, we can choose the parameters accordingly, and use the inclusion \eqref{AE},
 to obtain
 $$
\mathcal{A}_{q}^{\alpha} =
\,\big(\mathcal{A}_{\rho}^{\alpha_{0}},\mathcal{A}_{\rho}^{\alpha_{1}}\big)_{1/2,\,
q} \hookrightarrow
\big(\ell_{k^{\alpha_{0}}\eta(k)}^{q},\ell_{k^{\alpha_{1}}\eta(k)}^{q}\big)_{{1}/{2},\,
q} \,=\, \ell_{k^{\alpha}\eta(k)}^{q}(\mathcal{B},\mathbb{B}).
 $$

\sline ``$3\Rightarrow 4$''. This is trivial since
$\mathscr{G}_{q}^{\alpha}\hookrightarrow
\mathcal{A}_{q}^{\alpha}$.

\sline ``$4\Rightarrow 2$''. This is trivial from (\ref{4.1}).
\end{proof}

\begin{Remark} \label{RemProofThBernIne3}
{\rm Observe that $3\Rightarrow4\Rightarrow2\Leftrightarrow1$ hold
with the weaker assumption $\{k^\al\eta(k)\}\in\SW$.}
\end{Remark}

\begin{corollary}\label{CorProofThBernIne3} {\bf Optimal inclusions of \boldmath{$\cA_{q}^{\alpha}$} into \boldmath{$\ell^{q}_{\om}$}.}

\sline Let $\mathcal B$ be a (normalized) unconditional basis in
$\SB$. Fix $\alpha
>0$ and $q\in(0,\infty]$.

\sline (a) If  $h_\ell(N)$ is
doubling then $\mathcal{A}_{q}^{\alpha}(\mathcal{B},\mathbb{B})
\hookrightarrow
\ell_{k^{\alpha}h_\ell(k)}^{\,q}(\mathcal{B},\mathbb{B})$.

\sline (b) If for some $\om\in\SW$ we have $\mathcal{A}_{q}^{\alpha}(\mathcal{B},\mathbb{B})
\hookrightarrow
\ell_{\om}^{\,q}(\mathcal{B},\mathbb{B})$ then necessarily $\om(k)\lesssim k^\al h_\ell(k)$,
and hence $\ell_{k^{\alpha}h_\ell(k)}^{\,q}\hookrightarrow\ell_{\om}^{\,q}$.
\end{corollary}

\begin{proof}
Part (a) is an application of $1\Rightarrow3$ in the theorem with $\eta=h_\ell$
(which under the doubling assumption satisfies $\{k^\al h_\ell(k)\}\in\SW_+$ for all $\al>0$).
Part (b) is just a restatement of $3\Rightarrow1$ in the theorem,
setting $\eta(k)=\om(k)/k^\al$ and taking into account Remark \ref{RemProofThBernIne3}.
 \end{proof}


\section{Examples and Applications}  \label{Examples}

In this section we describe the democracy functions $h_\ell$ and
$h_r$ in various examples which can be found in the literature.
Inclusions for $\Aaq(\cB,\SB)$ and
$\mathscr{G}_{q}^{\alpha}(\cB,\SB)$ will be obtained inmediately
from the results of sections \ref{secRigtDemJackIneq} and
\ref{secLeftDemBernsIneq}. The most interesting case appears when
$\cB$ is a wavelet basis, and $\SB$ a function or distribution space
in $\SR^d$ which can be characterized by such basis (eg, the general
Besov  or Triebel-Lizorkin spaces, $B^\al_{p,q}$ and $F^{s}_{p,q}$,
and also rearrangement invariant spaces as the Orlicz and Lorentz
classes, $L^\Phi$ and $L^{p,q}$). Such characterizations provide a
description of each $\SB$ as a sequence space, so for simplicity we
shall work in this simpler setting, reminding in each case the
original function space framework.
\medskip

Let $\mathcal{D}=\cD(\SR^d)$ denote the family of all dyadic cubes
$Q$ in $\SR^d$, ie \[\cD=\big\{\,Q_{j,k} = 2^{-j}([0,1)^{d}+k)\mid
\mbox{\small$j\in\mathbb{Z},\;k\in\mathbb{Z}^{d}$}\,\big\}.\] We
shall consider sequences indexed by $\cD$,
$\bs=\{s_Q\}_{Q\in\cD}$, endowed with quasi-norms of the following
form \Be 
\Bigl\|\Big(\sum_{Q\in\cD}\big(\,|Q|^{\ga-\frac12}\,|s_Q|\,
{\chi_Q(\,\cdot\,)}
\,\big)^r\Big)^{1/r}\Bigr\|_{\SX}\;\;,
 \label{Xr} \Ee
where $0<r\leq\infty$, $\ga\in\SR$ and $\SX$ is a suitable quasi-Banach
function space in $\SR^d$, such as the ones we consider below. The
canonical basis $\cB_c=\{\be_Q\}_{Q\in\cD}$ is formed by the
sequences $\be_Q$ with entry 1 at $Q$ and 0 otherwise. In each of
the examples below, the greedy
algorithms and democracy functions are considered with respect to the normalized basis
$\cB=\big\{\be_Q/\|\be_Q\|_{\SB}\big\}$. Similarly, when stating the corresponding results
for the functional setting we shall write $\cW$ for the wavelet basis.


 \BE
{\boldmath{$\SX=L^p(\SR^d)$, {\small$0\!<\!p\!<\!\infty$}.}} In this
case, it is customary to consider the sequence spaces $\fpr$, $s\in
\mathbb R,\, 0 < r \leq \infty,$ with quasi-norms given by \[
\bigl\|\bs\bigr\|_{\fpr}\,:=\,
\Bigl\|\Big(\sum_{Q\in\cD}\big(|Q|^{-\frac s
d-\frac12}|s_Q|\,\chi_Q(\,\cdot\,)\,\big)^r\Big)^{1/r}\Bigr\|_{L^p(\SR^d)}\;\;.
\] It was proved in \cite{HJLY,GH,JM2} that, for all $s\in\SR$ and $0<r\leq\infty$,
\Be h_\ell(N;\fpr)\approx h_r(N;\fpr)\approx N^{1/p}\quad
\label{H_fpr}\Ee
 and \Be\Aaq(\fpr)=\ell^{\tau,q}(\fpr)=\Big\{\bs\mid \{s_Q\|e_Q\|_{\fpr}\}_Q\in\ell^{\tau,q}\Big\},\label{A_fpr}\Ee if
 $\frac1\tau=\al+\frac1p$,
 as asserted in Theorem \ref{KPGH}.

It is well-known that $\fpr$ coincides with the coefficient space
under a wavelet basis $\cW$ of the (homogeneous) Triebel-Lizorkin
space
 $\dot{F}_{p,r}^{s}(\mathbb{R}^{d}),$ defined in terms of Littlewood-Paley theory (see eg \cite{FJ,Meyer,Kyr}). In particular, under suitable decay and smoothness
on the wavelet family (so that it is an unconditional basis of the
involved spaces) the statement in \eqref{A_fpr} can be translated
into
  $$\mathcal{A}_{q}^{{\alpha}}(\mathcal{W},\dot{F}_{p,r}^{s}(\mathbb{R}^{d}))
 = \mathscr{G}_{q}^{{\alpha}}(\mathcal{W},\dot{F}_{p,r}^{s}(\mathbb{R}^{d}))
 = \dot{B}_{q,q}^{s+\alpha d}(\mathbb{R}^{d})
   $$
when $\frac{1}{q} = \frac{\alpha}{d} + \frac{1}{p}$. We refer to
\cite{HJLY,JM,DeV,GH} for details and further results.
\EE

{\bf Example 5.2.} \textbf{Weighted Lebesgue spaces
{\boldmath$\SX=L^p(w)$, $0<p<\infty$}}. For weights $w(x)$ in the
Muckenhoupt class $A_\infty(\SR^d)$, one can define sequence
spaces $\fpr(w)$ with the quasi-norm
\[ \bigl\|\bs\bigr\|_{\fpr(w)}\,:=\,
\Bigl\|\Big(\sum_{Q\in\cD}\big(|Q|^{-\frac s
d-\frac12}|s_Q|\,\chi_Q(\,\cdot\,)\,\big)^r\Big)^{1/r}\Bigr\|_{L^p(\SR^d,w)}\;\;.
\]
Similar computations as in the previous case in this more general
situation will also lead to the identities in \eqref{H_fpr} and
\eqref{A_fpr}, with $\fpr$ replaced by $\fpr(w)$. We refer to
\cite{Nat, KP2} for details in some special cases.

When $\cW$ is a (sufficiently smooth) orthonormal wavelet basis
and $w$ is a weight in the Muckenhoupt  class $A_p(\SR^d)$,
$1<p<\infty$, then $\mathfrak{f}^0_{p,2}(w)$ becomes the
coefficient space of the weighted Lebesgue space $L^p(w)$ (see eg
\cite{ABM}). One then obtains as special case
  $$ h_\ell(N;\mathcal{W},
L^{p}(w))\approx h_{r}(N;\mathcal{W},L^{p}(w))\approx
N^{\frac{1}{p}}\,.
  $$
  Moreover, if $\omega \in A_\tau(\mathbb R^d)\,,$
  $$ \mathcal{A}_{\tau}^{\alpha}(\mathcal{W},L^{p}(w))\approx {\mathscr G}_{\tau}^{\alpha}
 (\mathcal{W},L^{p}(w))\approx \dot{B}_{\tau,\tau}^{\alpha d}(w^{\tau/p}),\quad
 \mbox{if }\;\tfrac{1}{\tau} = {\alpha} + \tfrac{1}{p}\;,$$ where
  $\dot{B}_{\tau,q}^{\alpha}(w)$ denotes a weighted Besov space
  (see \cite{Nat} for details).\\

 {\bf Example 5.3.} \textbf{Orlicz spaces {\boldmath$\SX=L^\Phi(\SR^d)$}}.
 Following \cite{GHM}, we denote by $\fphi$ the sequence space with quasi-norm
\[
\|\bs\|_{\fphi}\,:=\,
\Bigl\|\Big(\sum_{Q\in\cD}\big(|s_Q|\,\tfrac{\chi_Q(\,\cdot\,)}{|Q|^{1/2}}\,\big)^2\Big)^{1/2}\Bigr\|_{L^\Phi(\SR^d)}\;\;,
\]
where $L^\Phi$ is an Orlicz space with non-trivial Boyd indices.
If we denote by $\varphi(t)=1/\Phi^{-1}(1/t)$, the fundamental
function of $L^\Phi$, then it is shown in \cite{GHM} that
 $$
  h_{\ell}(N;\fphi)
  \,\approx \,\inf_{s>0}\tfrac{\varphi(Ns)}{\varphi(s)}
  \quad\mand \quad h_r(N; \fphi)
\,  \approx \,\sup_{s>0}\tfrac{\varphi(Ns)}{\varphi(s)},
 $$
with the two expressions being equivalent iff $\varphi(t)=t^{1/p}$
(ie, iff $L^\Phi=L^p$). Thus, these are first examples of
non-democratic spaces, with a wide range of possibilities for the
democracy functions. The theorems in sections
\ref{secRigtDemJackIneq} and \ref{secLeftDemBernsIneq} recover the
embeddings obtained in \cite{GHM} for the approximation classes
$\Aaq(\fphi)$ and $\Gaq(\fphi)$ in terms of weighted discrete
Lorentz spaces. When using suitable wavelet bases, these lead to
corresponding inclusions for $\Aaq(\cW, L^\Phi)$ and $\Gaq(\cW,
L^\Phi)$, some of which can be expressed in terms of Besov
spaces of generalized smoothness (see \cite{GHM} for details).\\

  {\bf Example 5.4.}  \textbf{Lorentz spaces {\boldmath$\SX=L^{p,q}(\SR^d)$, {\small $0<p, q<\infty$}.}}
Consider sequence spaces $\lpr$ defined by the following
quasi-norms
\[
\|\bs\|_{\lpr}\,:=\,
\Bigl\|\Big(\sum_{Q\in\cD}\big(|s_Q|\,\tfrac{\chi_Q(\,\cdot\,)}{|Q|^{1/2}}\,\big)^2\Big)^{1/2}\Bigr\|_{L^{p,q}(\SR^d)}\;\;.
\]
Their democracy functions have been computed in \cite{HMN},
obtaining
$$
  h_\ell(N;\lpr)\,\approx \,N^{\frac{1}{\max(p,q)}}
  \quad\mand\quad h_{r}(N;\lpr)\,\approx\,
  N^{\frac{1}{\min(p,q)}}\,.
 $$
These imply corresponding inclusions for the classes
$\cA^\al_s(\lpr)$ and $\mathscr G^\al_s(\lpr)$ in terms of discrete
Lorentz spaces $\ell^{\tau,s}$ (as described in the theorems of
sections \ref{secRigtDemJackIneq} and \ref{secLeftDemBernsIneq}).
The spaces $\lpr$ characterize, via wavelets, the usual Lorentz
spaces
  $L^{p,q}(\mathbb{R}^{d})$ when {\small $1< p<\infty$ and $1\leq q<\infty$}  (\cite{Soardi}).
Hence inclusions for $\cA^\al_s(\cW, L^{p,q})$ and $\mathscr
G^\al_s(\cW, L^{p,q})$ can be obtained using standard Besov spaces.
\\

 {\bf Example 5.5.}  \textbf{Hyperbolic wavelets.}
For $0<p<\infty$, consider now the sequence space
\[
\|\bs\|_{\hpr}\,:=\,
\Bigl\|\Big(\sum_{R}\big(|s_R|\,\tfrac{\chi_R(\,\cdot\,)}{|R|^{1/2}}\,\big)^2\Big)^{1/2}\Bigr\|_{L^{p}(\SR^d)}\;\;.
\]
where $R$ runs over the family of all dyadic rectangles of
$\SR^d$, that is $R=I_1\times\ldots\times I_d$, with
$I_i\in\cD(\SR)$, $i=1,\ldots,d$. This gives another example of
non-democratic basis. In fact,
 the
  following result is proved in \cite[Proposition 11]{Wo} (see also \cite{Tem2}):

\sline (a) If $ 0<p\leq 2,$
 $$
   h_\ell(N;\hpr)\, \approx\,
 N^{1/p} (\log N)^{(\frac12 - \frac1p)(d-1)}\quad\mand\quad
 h_r(N;\hpr)\,
 \approx \,N^{1/p}.
 $$

\sline (b) If $2 \leq p < \infty$,
 $$
   h_\ell(N;\hpr) \,\approx\,
 N^{1/p} \quad\mand\quad  h_r(N;\hpr)\,
 \approx \,N^{1/p} (\log N)^{(\frac12 - \frac1p)(d-1)}\,.
 $$
If $\mathcal{H}_d$ denotes the multidimensional (hyperbolic) Haar
basis, then $\hpr$ becomes the coefficient space of the usual
$L^p(\SR^d)$ if {\small $1<p<\infty$} (and the dyadic Hardy space
$H^p(\SR^d)$ if $0<p\leq1$). In this case, one obtains corresponding
inclusions for the classes $\Aaq(\mathcal{H}_d,L^p)$ and $\mathscr
G_q^\alpha(\mathcal{H}_d,L^p)$ (see also \cite[Thm 5.2]{KT2}), some
of which could possibly be expressed in terms of Besov spaces of
bounded mixed smoothness \cite{KT2,DKT}.
\\

{\bf Example 5.6.}  \textbf{Bounded mean oscillation.} Let $bmo$
denote the space of sequences $\bs=\{s_I\}_{I\in\cD}$ with
 \begin{equation} \label{bmo}
 \|{\bf s}\|_{bmo} = \sup_{I\in \mathcal D} \Big(\frac{1}{|I|}
 \sum_{J\subset I\,,J\in \mathcal D} |s_J|^2 |J|\Big)^{1/2} < \infty \,.
 \end{equation}
This sequence space gives the correct characterization of
$BMO(\mathbb R)$ for sufficiently smooth wavelet bases appropriately
normalized(see \cite{W82,FJ,HJLY}). Their democracy functions are
determined by
 \begin{equation} \label{democracybmo}
h_\ell(N; bmo)\approx 1\,, \quad h_r(N; bmo) \approx \big(\log
N\big)^{1/2}\,.
 \end{equation}
The first part of (\ref{democracybmo}) is easy to prove, and the
second follows, for instance, by an argument similar to the one
presented in the proof of \cite[Lemma 3]{Oswald}. Our results of
sections \ref{secRigtDemJackIneq} and \ref{secLeftDemBernsIneq}
give in this case the inclusions:
 \begin{equation} \label{inclusionsbmo}
 \ell^q_{k^\alpha\sqrt{\log k}}\hookrightarrow
 \mathscr G^\alpha_q (bmo) \hookrightarrow
 \mathcal A^\alpha_q (bmo) \hookrightarrow
 \ell^q_{k^\alpha}\,=\,\ell^{1/\al,q}\,.
 \end{equation}
However, this is not the best one can say for the approximation
classes $\Aaq$. A result proved in \cite{RT} (see also Proposition
11.6 in \cite{HJLY}) shows that one  actually has  \[ \mathcal
A^\alpha_q (bmo) \,= \,\mathcal A^\alpha_q
(\ell^\infty)\,=\,\ell^{1/\al,q},\] for all $\al>0$ and
$q\in(0,\infty]$. For $0<r<\infty$ one can define the space $bmo_r$
replacing the 2 by $r$ in (\ref{bmo}); it can then be shown that
$h_r(N; bmo_r) \approx \big(\log N\big)^{1/r}\,$ and $ \mathcal
A^\alpha_q (bmo_r) \,= \,\ell^{1/\al,q}.$

\section{Democracy Functions for $\mathcal{A}_{q}^{\alpha}(\mathcal{B},\mathbb{B})$ and
$\mathscr{G}_{q}^{\alpha}(\mathcal{B},\mathbb{B})$}

As usual, we fix a (normalized) unconditional basis  $\mathcal B =
\{e_{j}\}_{j=1}^\infty$ in $\SB$. In this section we compute the
democracy functions for  the spaces
$\mathcal{A}_{q}^{\alpha}(\mathcal{B},\mathbb{B})$ and
$\mathscr{G}_{q}^{\alpha}(\mathcal{B},\mathbb{B})$, in terms of
the democracy functions in the ambient space $\mathbb{B}$. To
distinguish among these notions we shall use, respectively, the
notations
$$h_\ell(N;\mathcal{A}_q^\al),\quad
h_\ell(N;\mathscr{G}_q^\al)\mand h_\ell(N;\SB),$$ and similarly
for $h_r$ (recall the definitions in section 2.5). Since we shall
use the embeddings in sections 3 and 4, observe first that \Be
 h_{\ell}(N;
\ell_{\om}^{q}(\mathcal{B},\mathbb{B}))\approx
h_r(N;\ell_{w}^{q}(\mathcal{B},\mathbb{B}))\approx \om(N),
\label{lemEquDemFuncLrDis} \Ee for all $\om\in\SW_+$ and
$0<q\leq\infty$. This is immediate from the definition of the
spaces $\ell^q_\om(\cB,\SB)$ and Lemma \ref{lemanorma}.

\begin{proposition}\label{propEquDemFunGreS}
Fix $\alpha>0$ and $\,0<q\leq\infty$. If
$h_\ell(\,\cdot\,;\mathbb{B})$ is doubling then \Benu \item[(a)]
$\Ds h_\ell(N;\mathscr{G}_{q}^{\alpha})\approx
N^{\alpha}h_\ell(N;\mathbb{B})$. \ \item[(b)] $ \Ds
h_{r}(N;\mathscr{G}_{q}^{\alpha})\approx
 N^{\alpha}h_{r}(N;\mathbb{B})$.\Eenu
 In particular, $\mathcal B$ is democratic in $\mathscr{G}_{q}^{\alpha}(\mathcal B, \mathbb
 B)$ if and
 only if $\mathcal B$ is democratic in $\mathbb B$.
 \end{proposition}

 \begin{proof}\label{proPropEquDemGreS}
 The inequalities ``{\small $\gtrsim$}'' in (a), and ``{\small $\lesssim$}'' in (b) follow
 immediately from the embeddings\[
\ell_{k^\alpha h_r(k)}^{q}(\mathcal{B};\mathbb{B}) \hookrightarrow
\mathscr{G}_{q}^{\alpha}(\mathcal{B},\mathbb{B})
 \hookrightarrow \ell_{k^\alpha
 h_\ell(k)}^{q}(\mathcal{B};\mathbb{B})
 \]and the remark in \eqref{lemEquDemFuncLrDis}.
Thus we must show the converse inequalities. To establish (a),
given $N = 1,2,3,\ldots$ choose $\Gamma$ with $|\Gamma| = N$ and
so that $\|1_{\Gamma}\|_{\mathbb{B}}\leq2 h_\ell(N;\mathbb{B})$.
  Then, using the trivial bound in \eqref{4.1} we obtain
\begin{eqnarray*}
h_\ell(N;\mathscr{G}_{q}^{\alpha})\,\leq\,\|1_{\Gamma}\|_{\mathscr{G}_{q}^{\alpha}}
\lesssim\, N^{\alpha} \|1_{\Gamma}\|_{\mathbb{B}}  \,\approx\,
N^{\alpha}h_\ell(N;\mathbb{B}).
\end{eqnarray*}

We now prove ``{\small $\gtrsim$}'' in (b). Given $N =1,2,\ldots$,
choose first $\Gamma$ with $|\Gamma| = N$ and
   $\|1_{\Gamma}\|_{\mathbb{B}}\geq
 \tfrac12  h_{r}(N;\mathbb{B})$, and then any $\Gamma'$
   disjoint with $\Ga$ with $|\Gamma'| = N$.
   Then
   \[
h_{r}(2N;\mathscr{G}_{q}^{\alpha}) \, \geq \,
\big\|1_{\Ga\cup\Ga'}\big\|_{\Gaq}\,\gtrsim
 N^\al \gamma_N (1_{\Ga\cup\Ga' };\mathbb B) \gtrsim
N^\al\,\big\|1_{\Ga}\big\|_{\SB} 
\,\approx\, N^\al h_r(N;\mathbb B).
   \]
The required bound then follows from the doubling property of
$h_r$.
\end{proof}

\begin{proposition}\label{proEquAproxS}
Fix $\alpha>0$ and $\,0<q\leq\infty$, and assume that
$h_\ell(\,\cdot\,;\mathbb{B})$ is doubling. Then \Benu \item[(a)]
$\Ds h_\ell(N;\Aaq)\approx N^{\alpha}h_\ell(N;\mathbb{B})$. \
\item[(b)] $ \Ds h_{r}(N;\Aaq)\lesssim
 N^{\alpha}h_{r}(N;\mathbb{B})$.\Eenu
  In particular, if
$\mathcal{B}$ is democratic in $\mathbb{B}$ then $\mathcal{B}$ is
democratic in $\mathcal{A}_{q}^{\alpha}(\mathcal{B},\mathbb{B}).$
\end{proposition}

\begin{proof}
As before,  ``{\small $\gtrsim$}'' in (a), and ``{\small
$\lesssim$}'' in (b) follow immediately from the embeddings\[
\ell_{k^\alpha h_r(k)}^{q}(\mathcal{B};\mathbb{B}) \hookrightarrow
\AaqB
 \hookrightarrow \ell_{k^\alpha
 h_\ell(k)}^{q}(\mathcal{B};\mathbb{B}).
 \]
The converse inequality in (a) follows from the previous
proposition and the trivial inclusion $\Gaq\hookrightarrow\Aaq$.
\end{proof}

As shown in Example 5.6, the converse to the last statement in
Proposition \ref{proEquAproxS} is not necessarily true. The space
$\SB=bmo$ is not democratic, but their approximation classes
$\Aaq(bmo)=\ell^{1/\al,\,q}$ are democratic. Moreover, this
example shows that the converse to the inequality in (b) does not
necessarily hold, since
$$h_r(N;\mathcal A^q_\al(bmo)) = N^\alpha\quad\mbox{ but }\quad N^\al h_r(N; bmo)
\approx N^\alpha (\log N)^{1/2}.$$

Nevertheless, we can give a sufficient condition for
$h_r(N;\Aaq)\approx N^\al h_r(N;\SB)$, which turns out to be
easily verifiable in all the other examples presented in $\S5$.
\\

{\bf P{\footnotesize\bf ROPERTY} (H)}.\label{PropertyH} We say that
$\cB$ satisfies the \textbf{Property (H)} if for each $n=1,2,3,... $
there exist $ \Gamma_n \subset\SN$, with $|\Gamma_n| = 2^n$,
satisfying the property
  \[ \|1_{\Gamma'}\|_{\mathbb{B}}\approx
h_{r}(2^{n-1};\SB),\quad \forall\;\Gamma'\subset
\Gamma_n\quad\mbox{with}\quad |\Gamma'| =2^{n-1}.\]

\begin{proposition}\label{propPoprieH}
Assume that $\mathcal{B}$ satisfies the Property (H). Then, for
all $\alpha > 0$ and  $0<q\leq\infty$
  $$
  h_{r}(N;\mathcal{A}_{q}^{\alpha})\approx N^{\alpha}h_{r}
  (N;\mathbb{B})
  $$
  \end{proposition}

  \begin{proof}\label{proofPropProprieH}
 We must show {\small ``$\gtrsim$''}, for which we argue as in the proof of Proposition \ref{propEquDemFunGreS}.
 Given $N=2^n$, select  $\Gamma_n$ as in the definition of Property (H). Then,
\[
h_{r}(N;\Aaq) \, \geq \, \big\|1_{\Ga_n}\big\|_{\Aaq}\,\gtrsim
N^\al\,\sigma_{N/2}(1_{\Ga_n}).\] Now, the property (H) (and the
remark in \eqref{ErrorNterApIneq}) give
\[\sigma_{N/2}(1_{\Ga_n})\,=\,
  \inf 
  \,\bigl\{\|1_{\Gamma'}\|_{\mathbb{B}}\mid
\mbox{\small$\Gamma'\subset\Gamma,\;|\Gamma'| = N/2$}
\bigr\}\,\approx \,h_{r}(N/2;\mathbb{B})\approx
\,h_{r}(N;\mathbb{B}).
\]
Combining these two facts 
the proposition follows for $N=2^n$. For general $N$ use the result
just proved and the doubling property of $h_r$.
\end{proof}

As an immediate consequence, the property (H) allows to remove the
possible logarithmic loss for the embedding $\ell^q_{k^\al
h_r(k)}(\cB,\SB)\hookrightarrow\Aaq(\cB,\SB)$ discussed in
Corollary \ref{CorProofThJacIne2}.

\begin{corollary} {\bf{More about optimality for inclusions into
\boldmath{$\cA_{q}^{\alpha}$}.}}

\noindent Assume that $(\SB,\cB)$ satisfies property (H). If for
some $\al>0$, $q\in(0,\infty]$ and $\om\in\SW_+$ we have
$\ell_{\om}^{q}(\mathcal{B},\mathbb{B})
\hookrightarrow\cA_{q}^{\alpha}(\mathcal{B},\mathbb{B})$, then
necessarily $\om(k)\gtrsim k^\al h_r(k)$, and therefore
$\ell^q_\om\hookrightarrow\ell^q_{k^\al h_r(k)}$.
\end{corollary}

The following examples show that Property (H) is often satisfied.\\

{\bf Example 6.1.} Wavelet bases in  Orlicz spaces $L^\Phi(\mathbb
R^d)$ satisfy the property (H). Indeed, recall from \cite[Thm
1.2]{GHM} (see also Example 5.3) that
 \begin{equation} \label{Exampleh2}
 h_r(N; L^\Phi) \approx\;\sup_{s>0}
\,{\varphi(Ns)}/{\varphi(s)}\,.
 \end{equation}
Moreover, any collection $\Gamma$ of $N$ pairwise disjoint dyadic
cubes \emph{with the same fixed size} $a > 0$ satisfies
 \begin{equation} \label{Exampleh1}
 \|1_\Gamma\|_{L^\Phi}\; \approx \,{\varphi(Na)}/{\varphi(a)}\,,
 \end{equation}
 (see eg
\cite[Lemma 3.1]{GHM}). Thus, for each $N=2^n$, we first select
$a_n=2^{j_n d}$ so that  $h_r(2^n;L^\Phi) \approx {\varphi(2^n
a_n)}/{\varphi(a_n)}\,$, and then we choose as $\Ga_n$ any
collection of $2^n$ pairwise disjoint cubes with constant size
$a_n$. Then, any subfamily $\Gamma' \subset \Gamma_n$ with
$|\Gamma'| = N/2$, satisfies
$$\|1_{\Gamma'}\|_{L^\Phi}\approx \;{\varphi((N/2)a_n)}/{\varphi(a_n)}
\approx
\;{\varphi(Na_n)}/{\varphi(a_n)}\,\approx\,h_r(N)\,\approx\,h_r(N/2),$$
by
\eqref{Exampleh1} and the doubling property of $\varphi$ and $h_r$. \\

{\bf Example 6.2.}  Wavelet bases in Lorentz spaces
$L^{p,q}(\mathbb R^d)$,  $1<p,q<\infty$. These also satisfy the
property (H). Indeed, it can be shown that any set $\Ga$
consisting of $N$ \emph{disjoint cubes of the same size} has
\[
\|{1}_{\Gamma}\|_{L^{p,q}}\,\approx\, N^{\frac{1}{p}}\;,
\]
while sets $\Dt$ consisting of $N$ disjoint cubes
 \emph{all having different sizes} satisfy \[
\|{1}_{\Dt}\|_{L^{p,q}}\,\approx\, N^{\frac{1}{q}}\;.
\]
(see \cite[(3.6) and (3.8)]{HMN}). Since $h_r(N)\approx
N^{1/(p\wedge q)}$, we can define the $\Ga_n$'s with sets of the
first type when $p\leq q$, and with sets of the second type when
$q<p$, to obtain in both cases a collection satisfying the
hypotheses of
 property (H).
\\

{\bf Example 6.3.} The hyperbolic Haar system in $L^p(\SR^d)$ from
Example 5.5 also satisfies property (H). In this case, again, any
set $\Ga$ consisting of $N$ disjoint rectangles has
\[
\|{1}_{\Gamma}\|_{L^{p}(\SR^d)}\,=\, N^{\frac{1}{p}}\;.
\]
On the other hand, if $\Dt_n$ denotes the set of all the dyadic
rectangles in the unit cube with fixed size $2^{-n}$, then \Be
\|{1}_{\Dt_n}\|_{L^{p}(\SR^d)}\,\approx\,
2^{n/p}\,n^{(d-1)/2}\,\approx\,|\Dt_n|^{1/p}(\log
|\Dt_n|)^{(d-1)(\tfrac12-\tfrac1p)}\;. \label{Dtp}\Ee Moreover, it
is not difficult to show that any $\Dt'\subset\Dt_n$ with
$|\Dt'|=|\Dt_n|/2$ also satisfies \eqref{Dtp} (with $\Dt_n$
replaced by $\Dt'$). Hence, combining these two cases and using
the description of $h_r(N)$ in Example 5.5, one easily establishes
the property (H).

\section{Counterexamples for the classes $\GaqB$}
\subsection{Conditions for
{\boldmath $\Gaq\not=\Aaq$} } \label{Noninclusions}

Recall from section 2.3 that
$\mathscr{G}_{q}^{\alpha}(\mathcal{B},\mathbb{B})\hookrightarrow
\mathcal{A}_{q}^{\alpha}(\mathcal{B},\mathbb{B})$, with equality
of the spaces when $\cB$ is a democratic basis. It is known that
there are some \emph{conditional} non-democratic bases for which
$\Gaq=\Aaq$ (see \cite[Remark 6.2]{GN}). For unconditional bases,
however, one could ask whether non-democracy
 necessarily implies that $\Gaq\not=\Aaq$. We do
not know how to prove such a general result, but we can show that
the inclusion $\Aaq\hookrightarrow\Gaq$ must fail whenever the gap
between $h_\ell(N)$ and $h_r(N)$ is at least logarithmic (and even
less than that). More precisely, we have the following.

\begin{proposition} \label{prop7.1}
Let $\mathcal{B}$ be an unconditional basis in $\SB$ and $\alpha
> 0$. Suppose that there exist integers $p_N \geq q_N\geq1$, $N= 1, 2, \dots $
 such that
 \begin{eqnarray} \label{eq7.1}
  \lim_{N\to \infty}
\frac{p_N}{q_N} = \infty \quad\mand\quad
\frac{h_r(q_N)}{h_\ell(p_N)} \gtrsim
\left(\frac{p_N}{q_N}\right)^\alpha\,.
 \end{eqnarray}
 Then the inclusion
$\mathcal{A}_{\tau}^{\alpha}(\mathcal{B},\mathbb{B})\hookrightarrow
\mathscr{G}_{\tau}^{\alpha}(\mathcal{B},\mathbb{B})$ does not hold
for any  $\tau\in(0,\infty].$
\end{proposition}

\begin{proof}
For each $N$, choose $\Ga_l,\Gamma_r \subset \mathbb N$ with
$|\Gamma_l| = p_N$, $|\Gamma_r| = q_N$, and such
 that
 \begin{eqnarray} \label{eq7.2}
\|1_{\Gamma_l}\|_{\mathbb B} \leq 2 h_\ell(p_N),\,\quad
\|1_{\Gamma_r}\|_{\mathbb B} \geq\, \tfrac12 \,h_r(q_N)\,.
 \end{eqnarray}
Set $x_N = \bone_{\Gamma_r} + 2\cdot \,\bone_{\Gamma_l - \Gamma
_l\cap \Gamma_r}\,.$ Since $\#(\Gamma_l - \Gamma _l\cap
\Gamma_r)\geq
p_N-q_N$, when $k \in [1, p_N - q_N]$ we have 
 $$
 \|x_N - G_k(x_N)\|_{\mathbb B} \geq \|
 1_{\Gamma_r}\|_{\mathbb B} \geq \,\tfrac12 \,h_r(q_N)\,.
 $$
Therefore, using $p_N - q_N > p_N/2$ (since ${p_N}/{q_N}>2$ for
$N$ large), we obtain that
 \begin{eqnarray} \label{eq7.4}
 \|x_N\|_{\mathscr{G}_\tau^\alpha (\mathcal B, \mathbb B)} \geq
\,\tfrac12\, \Big[\sum_{k=1}^{{p_N}/{2}} \big(k^\alpha
h_r(q_N)\big)^\tau\, \tfrac{1}{k}\Big]^{\frac{1}{\tau}} \gtrsim\,
h_r(q_N) \,p_N^\alpha\;.
 \end{eqnarray}
On the other hand, we can estimate the norm of $x_N$ as follows:
 \begin{eqnarray} \label{eq7.5}
 \| x_N \|_{\mathbb B} \lesssim \| 1_{\Gamma_r}\|_{\mathbb
 B} +  \|1_{\Gamma_l - \Gamma _l\cap
\Gamma_r}\|_{\mathbb B} \leq
 h_r(q_N) + 2 h_\ell(p_N)
 \lesssim h_r(q_N)
 \end{eqnarray}
where the last inequality is true for $N$ large due to
(\ref{eq7.1}). Thus
 \begin{eqnarray} \label{eq7.6}
 \sigma_k (x_N) \leq \| x_N \|_{\mathbb B}  \lesssim h_r(q_N)\,.
 \end{eqnarray}
Next, if $k \geq q_N$, by (\ref{eq7.2})
 \begin{eqnarray} \label{eq7.7}
 \sigma_k (x_N) \leq 2\|1_{\Gamma_l - \Gamma _l\cap
\Gamma_r}\|_{\mathbb B}
 \leq 2\|1_{\Gamma_l}\|_{\mathbb B}  \lesssim h_\ell(p_N)\,.
 \end{eqnarray}
Combining (\ref{eq7.5}), (\ref{eq7.6}), and (\ref{eq7.7}) we see
that
 \begin{eqnarray} \label{eq7.8}
 \|x_N\|_{\mathcal{A}_\tau^\alpha (\mathcal B, \mathbb B)}
 &\lesssim&
 h_r(q_N) + \Big[\sum_{k=1}^{q_N-1} \big(k^\alpha h_r(q_N)\big)^\tau \tfrac{1}{k} +
 \sum_{k=q_N}^{p_N + q_N} \big(k^\alpha h_\ell(p_N)\big)^\tau \tfrac{1}{k}
 \Big]^{\frac{1}{\tau}} \nonumber\\
 &\lesssim & h_r(q_N) + \big[ h_r(q_N)^\tau (q_N)^{\alpha\tau} +
 h_\ell(p_N)^{\tau} (p_N)^{\alpha\tau}\big]^{\frac{1}{\tau}}\nonumber\\
 &\lesssim&  h_r(q_N) + h_r(q_N) (q_N)^\alpha \lesssim h_r(q_N) (q_N)^\alpha
 \end{eqnarray}
where in the second inequality we have used the elementary fact
$\sum_{k=a}^{a+b} k^{\ga-1}\lesssim b^\ga$ if $b\geq a$, and the
third inequality is due to (\ref{eq7.1}). Therefore, from
(\ref{eq7.4}) and (\ref{eq7.8}) we deduce
 $$
 \frac{\|x_N\|_{\mathscr{G}_\tau^\alpha }}
 {\|x_N\|_{\mathcal{A}_\tau^\alpha}}\,
 \gtrsim \,\frac{h_r(q_N) (p_N)^\alpha}{h_r(q_N) (q_N)^\alpha} \,= \,\Big(\frac{p_N}{q_N}\Big)^\alpha
 \longrightarrow \infty
 $$
as $N \rightarrow \infty$. This shows the desired result.
\end{proof}

\begin{corollary}\label{cor7.2}
Let $\mathcal{B}$ be an unconditional basis such that
$h_\ell(N)\lesssim N^{\beta_0}$ and $h_{r}(N)\gtrsim N^{\beta_1}$,
for some $ \beta_1>\beta_0\geq0.$ Then,
$\Gaq\not=\mathcal{A}_{\tau}^{\alpha}$, for all $\alpha
> 0$ and all $\tau\in(0,\infty].$
\end{corollary}

\begin{proof}
Choose $ r, s \in \mathbb N\,,$ such that $\frac{\alpha +
\beta_0}{\alpha + \beta_1} <  \frac{r}{s} < 1.$ Take $p_N = N^s$ and
$q_N = N^r$. Then,
 $
 \lim_{N\to\infty} \frac{p_N}{q_N} = \lim_{N\to\infty} N^{s-r} =
 \infty
 $
 and
 $$
 \frac{h_r(q_N)}{h_\ell(p_N)} \gtrsim \frac{N^{r \beta_1}}{N^{s \beta_0}}
 > N^{\alpha(s-r)} =
 \Big(\frac{N^s}{N^r}\Big)^\alpha =
 \Big(\frac{p_N}{q_N}\Big)^\alpha\,,
 $$
which proves (\ref{eq7.1}) in this case, so that we can apply
Proposition \ref{prop7.1}.
\end{proof}

\begin{corollary}\label{cor7.3}
Let $\mathcal{B}$ be an unconditional basis such that for some
$\beta \geq 0 $ and $\gamma > 0$ we have either

\smallskip {\rm (i)} $\;h_{r}(N)\gtrsim N^{\beta} (\log N)^\gamma$
and $h_\ell(N)\lesssim N^{\beta},$
 or

\smallskip {\rm (ii)} $\;h_{r}(N)\gtrsim N^{\beta} $ and
$h_\ell(N)\lesssim N^{\beta}(\log N)^{-\gamma}.$

\sline Then, $\Gaq\not=\Aaq$ for all $\alpha > 0$ and all
$\tau\in(0,\infty].$
\end{corollary}

\begin{proof}
$i)$ Choose $a,b \in \mathbb N$ such that $0 < \frac{a}{b} <
\frac{\gamma}{\alpha + \beta}.$   Let $p_N = N^a 2^{N^b}$ and $ q_N
= 2^{N^b}$. Then,
 $
 \lim_{N\to\infty} \frac{p_N}{q_N} = \lim_{N\to\infty} N^{a} =
 \infty
 $
 and
 $$
 \frac{h_r(q_N)}{h_\ell(p_N)} \,\gtrsim\,\frac{(2^{N^b})^\beta (\log
2^{N^b})^\gamma}{N^{a\beta} (2^{N^b})^\beta}
 \thickapprox \frac{N^{b\gamma}}{N^{a\beta}} = N^{b\gamma - a\beta} >
 N^{a\alpha} = \Big(\frac{p_N}{q_N}\Big)^\alpha
 $$
which proves (\ref{eq7.1}) in this case, so that we can apply
Proposition \ref{prop7.1} to conclude the result. The proof of $
ii)$ is similar with the same choice of $p_N$ and $q_N.$
\end{proof}


\subsection{Non linearity of $\GaqB$}

We conclude by showing with simple examples that $\Gaq (\mathcal
B, \mathbb B)$ may not even be a linear space when the basis $\cB$
is not democratic.

Let $\mathbb B = \ell^p \oplus_{\ell^1} \ell^q$, $ 0 < q < p <
\infty$; that is, $\mathbb B$ consists of pairs $(a,b)\in
\ell^p\times\ell^q$, endowed with the quasi-norm $\|a\|_{\ell^p} +
\|b\|_{\ell^q}\,.$ We consider the canonical basis in $\SB$.

\smallskip

Now, set $\beta = \alpha +\frac{1}{p}\,$ and $x = \{(k^{-\beta},
0)\}_{k\in\mathbb N}\in\SB$. For $N=1,2,3, \dots$ we have
 $$ \gamma_N (x) = \Big(\sum_{k>N}
 \frac{1}{k^{\beta p}}\Big)^{1/p} \approx \Big(\frac{1}{N^{\beta
 p -1}} \Big)^{1/p} \,=\, N^{-\alpha}\,.
 $$
This shows that $x\in \mathscr{G}_\infty^\alpha (\mathcal B, \mathbb
B)$. Similarly, if we let $\gamma = \alpha +\frac{1}{q}$, then  $y =
\{(0 , j^{-\gamma})\}_{j\in\mathbb N}$ belongs to
$\mathscr{G}_\infty^\alpha$.
 We will show, however, that $x+y\not\in\mathscr{G}^\al_\infty$.
 In fact, we will find a subsequence $N_J$ of natural numbers so
 that
 \begin{equation}  \label{sum}
 \gamma_{N_J} (x+y) \approx \frac{1}{N_J^{\alpha \beta/\gamma}}
 \end{equation}
(notice that $\beta<\ga$ since we chose $q<p$). To prove
(\ref{sum}) let $A_1=\{1\}$ and
 $$
  A_j = \Big\{k \in \mathbb N : \frac{1}{j^\gamma} \leq
  \frac{1}{k^\beta} < \frac{1}{(j-1)^\gamma}\Big\}\,, \quad j= 2,3
  ,\dots
 $$
The number of elements in $A_j$ is
 \begin{equation} \label{number}
 |A_j| \,\approx\, j^{\gamma/\beta} - (j-1)^{\gamma/\beta} \approx
 j^{\frac{\ga}{\beta}-1}\,, \quad j=1,2,3, \dots
 \end{equation}
For $J =2,3,4, \dots$ let $N_J = \sum_{j=1}^J |A_j| + J$. From
(\ref{number}) we obtain
 $$
 N_J \approx \sum_{j=1}^J j^{\frac{\gamma}{\beta}-1} + J \approx
 J^\frac{\gamma}{\beta} + J \approx J^\frac{\gamma}{\beta}\,,
 $$
since $\gamma > \beta \,.$ Thus,
 \begin{eqnarray*}
 \gamma_{N_J}(x+y) &\approx & \Big( \sum_{k>
 J^{\frac\gamma\beta}} {k^{-\beta p}}\Big)^{1/p} + \Big( \sum_{j>
 J} {j^{-\gamma q}}\Big)^{1/q} \,\approx\,
 \left[(J^{\gamma/\beta})^{-\beta p +1}\right]^{1/p} +
 \left[J^{-\gamma q + 1}\right]^{1/q}\\
&  =& J^{-\al\ga/\beta}+J^{-\al}
 \;\approx  J^{-\alpha}\;
\; \approx \;(N_J)^{-\alpha \beta /\gamma}\,,
 \end{eqnarray*}
proving (\ref{sum}).

 A simple modification of the above construction can be used
to show that the set $\mathscr{G}_s^\alpha (\mathcal B, \mathbb
B)$ is not linear, for any $\alpha > 0$ and any $s\in(0,\infty)$.

\end{document}